\documentclass[12pt,twoside]{amsart}\usepackage{amssymb}
\usepackage{enumerate}

\newtheorem{thmspec}{\relax}

\newtheorem{theorem}{Theorem}[section]

\newtheorem{prop}[theorem]{Proposition}
\newtheorem{defi}[theorem]{Definition}

\theoremstyle{definition}

\theoremstyle{remark}

\numberwithin{equation}{section}
\def \cal{\mathcal}
\def \Bbb{\mathbb}

\def\onto{{\kern3pt\to\kern-8pt\to\kern3pt}}

\def\<{\langle}
\def\>{\rangle}
\def\|{{\ |\ }}

\def\onto{\twoheadrightarrow}

\def\-{\underline}

\def\dist{\operatorname{dist}}

\def\mes{\operatorname{mes}}

\def\arg{\operatorname{arg}}

\def\End{\operatorname{End}}
\def\Alim{\mathcal{A}-\lim}
\def\Alimsup{\mathcal{A}-\limsup}

\def\R{\Bbb R}
\def\C{\Bbb C}
\def\P{\Bbb P}

\def\S{\Bbb S}
\def\X{\Bbb X}

\def\<{\langle}
\def\>{\rangle}

\catcode`\@=11
\def\serieslogo@{\relax}
\def\@setcopyright{\relax}
\catcode`\@=12

\title[Recent developments]
{Recent developments in the theory of   separately holomorphic mappings}

\begin{document}

\author{Vi{\^e}t-Anh  Nguy\^en}
\address{Vi{\^e}t-Anh  Nguy\^en\\
Vietnamese Academy  of Science  and  Technology\\
Institute of Mathematics,
Department  of Analysis\\
18  Hoang Quoc  Viet  Road, Cau Giay  District\\
10307 Hanoi, Vietnam}
\email{nvanh@math.ac.vn}

\address{{\sc Current  address:}
 School of Mathematics\\
Korea Institute  for Advanced Study\\
 207-43Cheongryangni-2dong,
Dongdaemun-gu\\
Seoul 130-722, Korea} \email{vietanh@kias.re.kr}

\date{}


\begin{abstract}
We describe  a part of the recent developments
in the theory of   separately holomorphic mappings between complex analytic spaces.
Our description focuses on  works  using the technique of holomorphic discs.
\end{abstract}
\maketitle

\section{Introduction} \label{sec_introduction}

 In this exposition all complex manifolds are assumed to be  of  finite  dimension  and  countable at infinity, 
 and all
  complex analytic spaces  are assumed to be reduced,
  irreducible  and  countable at infinity.
For a subset $S$ of a topological space $\mathcal{M},$  $\overline{S}$   denotes the closure of $S$
in $\mathcal{M}.$   For  two  complex analytic spaces (resp. two topological spaces)  $D$  and  $Z,$ $\mathcal{O}(D,Z)$ (resp. $\mathcal{C}(D,Z)$) denotes the set of
all holomorphic (resp. continuous) mappings from  $D$ to $Z$.

\smallskip

The main purpose of this work  is to describe the recent developments
around  the following two problems.

\medskip

\noindent{\bf  PROBLEM 1.} {\it
 Let $X,$ $Y$ be two complex manifolds, let  $D$ (resp. $G$) be an open subset of $X$  (resp.  $Y$), let
  $A$ (resp. $B$) be a  subset of  $\overline{D}$ (resp.
  $\overline{G}$)  and let  $Z$ be
  a complex analytic  space.
 Define the {\bf cross}
\begin{equation*}
 W:=  \big((D\cup A)\times B\big) \bigcup\big(A\times (G\cup B)\big)  .
 \end{equation*}
   We  want to determine  the {\bf envelope of holomorphy} of the cross $W, $ that is,
   an ``optimal" open subset  of $ X\times Y,$ denoted by  $\widehat{\widetilde{W}},$ which is characterized by the following properties:

   For every mapping $f:W\longrightarrow Z$ that satisfies, in essence, the following condition:

$$
\begin{aligned}
 f(a,\cdot)\in&\mathcal{C}(G\cup B,Z)\cap\mathcal{O}(G,Z),\quad a\in A, \\
 f(\cdot,b)\in&\mathcal{C}(D\cup A,Z)\cap\mathcal{O}(D,Z),\quad b\in B,
\end{aligned}
$$

 there exists an  $\hat{f}\in\mathcal{O}(\widehat{\widetilde{W}},Z)$  such that for every $(\zeta,\eta)\in W,$
  $\hat{f}(z,w)$ tends to $f(\zeta,\eta)$ as $(z,w)\in\widehat{\widetilde{W}}$ tends, in some sense, to $(\zeta,\eta).$
  }

  \medskip

  The second problem generalizes the first one   to the case where we   add a set of singularities $M$  to the cross.
In order to understand this problem we need to introduce  some more notation and terminology.
Let $X,$  $Y,$ $D,$  $G,$ $A,$  $B$  and $Z$ and  $W$ be  as in PROBLEM 1
 and  let $M\subset W.$
The set $M_a:=\{w\in G:\ (a,w)\in M\}$,  $a\in A$, is  called
 \textit{the vertical fiber of  $M$ over  $a$} (resp. the set
$M^{b}:=\{z\in D:\ (z,b)\in M\}$, $b\in B$, is  called \textit{the
horizontal fiber  of $M$  over $b$}). We say that  $M$ possesses a
certain property    {\it in fibers over  $A$   (resp. $B$)} if all
vertical fibers $M_a$, $a\in A$,  (resp. all horizontal fibers
$M^b$, $b\in B$) possess  this property.

\medskip

\noindent {\bf  PROBLEM 2.} {\it
 Under  the  above hypotheses and notation  let  $\widehat{\widetilde{W}}$ be  the envelope of holomorphy
 of $W$  given by {\rm PROBLEM 1.}
 For  every  subset $M\subset W$   which is  relatively closed and  locally  pluripolar   (resp. thin)
\footnote{    The notion  of local pluripolarity (resp. thinness)  will be  defined  in Subsection
  \ref{Subsection_approach_regions} (resp. Section  \ref{sec_PROBLEM_2}) below.}
in fibers  over $A$  and $B$ ($M=\varnothing$ is  allowed)  we want to know  if there  exists an ``optimal"   set of singularities $\widehat{M} \subset \widehat{\widetilde{W}},$
which is relatively closed   locally  pluripolar   (resp. relatively closed analytic)
 and
which
 is characterized by the following property:

 For every mapping $f:W\setminus M\longrightarrow Z$ that satisfies, in essence, the following condition:

$$
\begin{aligned}
 f(a,\cdot)\in&\mathcal{C}((G\cup B)\setminus M_a,Z)\cap\mathcal{O}(G\setminus M_a,Z),\quad a\in A, \\
 f(\cdot,b)\in&\mathcal{C}((D\cup A)\setminus M^b,Z)\cap\mathcal{O}(D\setminus M^b,Z),\quad b\in B,
\end{aligned}
$$
 there exists an  $\hat{f}\in\mathcal{O}(\widehat{\widetilde{W}}\setminus \widehat{M},Z)$  such that
for all $(\zeta,\eta)\in W\setminus M$, $\hat{f}(z,w)$ tends to
$f(\zeta,\eta)$ as $(z,w)\in\widehat{\widetilde{W}}\setminus\widehat{M}$ tends, in some
sense, to $(\zeta,\eta)$. }

  \medskip

  The motivation  for PROBLEM 2  will be explained  in  Section \ref{sec_history}
  and  \ref{sec_PROBLEM_2} below.
  These  problems  play a fundamental role  in   the theory of separately holomorphic (resp. meromorphic)
 mappings, and they have been intensively studied during the last decades.
  There are two recent surveys by  Nguy\^en Thanh V\^an (see \cite{ng}) and by Peter  Pflug (see \cite{pf}) which
      summarize the  historical developments  up to  2001 of  PROBLEM 1 and 2 under
  the hypotheses that

  \smallskip

  {\bf $A\subset D$ and $B\subset G$ and $X,\ Y$ are Stein manifolds and  $Z$  is  a complex analytic space which possesses the Hartogs extension property \footnote{  This notion will be  defined in  Subsection  \ref{HEP} below.}.}

  \smallskip

      Both survey articles    give interesting insights
  and  suggest  new research trends in this subject.
  Our exposition  may be considered as a continuation to  the above  works.
  Namely, we describe  a part of the recent developments using the technique of holomorphic discs. This  will permit  us
 to obtain  partial (but
reasonable) solutions to  PROBLEM 1 and 2 in the case where  $Z$  is  a complex analytic space which possesses the Hartogs extension property.

\smallskip

  We close  the introduction with a brief outline of the paper to follow.

\smallskip

 In Section 2 we  describe   briefly    the  historical developments of PROBLEM 1 and 2.

In Section 3 we
 provide the  framework for an  exact  formulation of both
problems and for  their solution.

 The technique of holomorphic discs and related results
 are described
in Section 4.

In Section 5 we present some ideas of our new approach to the theory of separate holomorphy.
More precisely,   we apply  the results of Section 4 in order to complete  PROBLEM 1
in a  special case.

Section 6 is devoted to various  partial  results of PROBLEM 1.

Some approaches to  PROBLEM 1 and 2 are given in Section 7 and 8 respectively.
In fact,  Section 6 and 8 are  obtained  in collaboration  with Pflug  (see  \cite{pn1,pn2,pn3,pn4,pn5,np1}).

 Various applications of our solutions  are given in Section 9.

 Section 10 concludes the article
 with some  remarks and  open questions.

 \smallskip

\indent{\it{\bf Acknowledgment.}}
 The paper was written while  the  author was visiting  the  Abdus Salam International Centre
 for Theoretical Physics
in Trieste and the Korea Institute for Advanced Study in Seoul. He wishes to express his gratitude to these organizations.

\section{History} \label{sec_history}

Now we  recall    briefly  the main   developments around PROBLEM 1 and 2.
 All the results  obtained  so far  may be divided into two directions.
The first direction  investigates  the results in the  ``interior" context:
 $A\subset D$ and  $B\subset G,$  while
the second one explores the ``boundary"  context:  $A\subset\partial  D$  and  $B\subset\partial G.$

 \smallskip

The first  fundamental result in the field of separate holomorphy is
   the well-known Hartogs extension theorem for separately holomorphic functions
(see \cite{ha}). In the language  of  PROBLEM 1
 the following case:   $X=\C^n,\  Y=\C^m,\  A=D,\  B=G,  \ Z=\C$  has been solved,  and the
result is  $\widehat{\widetilde{W}}=D\times G.$
In particular, this theorem may be considered as the first main result in the first direction.
In 1912  Bernstein obtained, in his famous  article \cite{ber}, a positive solution to PROBLEM 1
for certain cases  where  $A\subset D,$  $B\subset G,$ $X= Y=\C$ and  $Z=\C.$

\smallskip

The next important development came about very much later.  In   1969--1970
 Siciak   established some significant  generalizations  of
the Hartogs
extension theorem  (see \cite{si1,si2}).
   In fact,    Siciak's   formulation of these generalizations gives rise to
  PROBLEM 1:    to determine the envelope of holomorphy for  separately
holomorphic functions defined on some {\it cross sets}  $W.$ The theorems
obtained under this formulation are often called {\it cross theorems.}
 Using  the so-called {\it relative
extremal function} (see Section \ref{sec_new_formulations} below),  Siciak completed  PROBLEM 1 for  the case where
$A\subset D,$  $B\subset G,$ $X= Y=\C$ and  $Z=\C.$

\smallskip

The next deep steps were initiated by Zahariuta in 1976 (see \cite{za}) when he started to use
the method of common bases of Hilbert spaces. This original approach permitted
him to obtain new  cross theorems for some cases where $A\subset D,\  B\subset G$ and $D=X,$  $G=Y$ are  Stein manifolds.
As a consequence, he was able to generalize
the result of Siciak in higher dimensions.

\smallskip

Later, Nguy\^en  Thanh V\^an and Zeriahi (see \cite{nz1,nz2,nz3}) developed the method of doubly orthogonal bases of
Bergman type in order to generalize the result of
Zahariuta. This  is a significantly simpler and    more constructive version
of Zahariuta's original  method.  Nguy\^en  Thanh V\^an and Zeriahi have recently achieved an elegant
 improvement  of their method  (see \cite{ng2}, \cite{zer}).

\smallskip

Using Siciak's  method, Shiffman  (see \cite{sh1}) was the first to generalize some
results of Siciak to separately holomorphic mappings with values in a complex analytic space  $Z.$
Shiffman's result (see \cite{sh2}) shows that the natural ``target spaces" for obtaining
satisfactory generalizations of cross theorems are the ones which possess {\it the Hartogs extension property}
(see Subsection  \ref{HEP} below for more explanations).

\smallskip

In 2001 Alehyane and Zeriahi   solved   PROBLEM 1  for the case where $A\subset D,$  $B\subset G$
 and $X,\ Y$  are Stein manifolds
and $Z$ is a complex analytic space which possesses the Hartogs extension property. The envelope of holomorphy
$\widehat{\widetilde{W}}$ is then given
 by
 \begin{equation*}
 \widehat{\widetilde{W}}
 :=\left\lbrace (z,w)\in D\times G :\  \widetilde{\omega}(z,A,D)+  \widetilde{\omega}(w,B ,G)<1
\right\rbrace,
\end{equation*}
 where $ \widetilde{\omega}(\cdot,A,D)$ and $ \widetilde{\omega}(\cdot,B ,G)$ are the plurisubharmonic measures,
 which  are generalizations of Siciak's relative extremal function
 (see Section  \ref{sec_new_formulations} below for this notion).
 This is the most general result to PROBLEM 1 under the hypothesis $A\subset D,$  $B\subset G.$
 More precisely,

 \medskip

\noindent {\bf Theorem 1   (Alehyane--Zeriahi \cite{az}).}
{\it
  Let $X,\ Y$ be   Stein manifolds, and $D\subset X,$  $G\subset Y$      domains,
  and  $ A\subset
D,$  $B\subset G$  nonpluripolar   subsets.
 Let $Z$ be a complex  analytic space possessing the
  Hartogs
  extension property.
Then for every mapping $f$ as in the hypotheses of {\rm PROBLEM 1,} there is a unique mapping
$\hat{f}\in\mathcal{O}(\widehat{\widetilde{W}},Z)$
such that $\hat{f}=f$ on $W\cap \widehat{\widetilde{W}}.$
}

\smallskip

In fact, Theorem 1 is  still valid for  $N$-fold crosses $W$ $(N\geq 2)$. For the notion of an $N$-fold cross
see, for example, \cite{pf} or \cite{nv1}.

 PROBLEM 2 has started with  a paper by \"{O}ktem in $1998$ (see \cite{ok1,ok2}) investigating the range problem in
Mathematical Tomography.
The reader will find in Section  \ref{sec_PROBLEM_2} below   a
 concise description of the range problem  and its relations to the theory of separate holomorphy.
  On the other hand,  Henkin and  Shananin gave, in an earlier work \cite{hs}, some applications
of Bernstein's result  \cite{ber} to Mathematical Tomography.
Here is the most general result in this direction. In fact, we state it in a somewhat simplified from.

\medskip

\noindent {\bf Theorem 2  (Jarnicki--Pflug  \cite{jp3,jp5}).}
  {\it
  Let $X$  and  $Y$ be Riemann--Stein domains, 
  let $D\subset X,$ $ G\subset Y$ be two subdomains, let
  $A\subset D$  and $B\subset G$ be nonpluripolar  subsets.
  Suppose in addition that $\mathcal{O}(D,\C)$ (resp. $\mathcal{O}(G,\C)$) separates points in $D$  (resp. in $G$)
  \footnote{ We say that   {\it $\mathcal{O}(D,\C)$  separates points in $D$} if for all points $z_1,$  $z_2$ with
  $z_1\not=z_2,$ there exists  $f\in  \mathcal{O}(D,\C)$ such that $f(z_1)\not=f(z_2).$}.
  Let  $M\subset W$ be a relatively closed subset which is    pluripolar (resp. thin)
  in fibers over $ A$ and  $ B.$

   Then there exists a relatively closed pluripolar  set  (resp.  relatively closed  analytic set) $\widehat{M}\subset \widehat{\widetilde{W}}$
   such that:
   \begin{itemize}
   \item[$\bullet$]
   $\widehat{M}\cap W\cap  \widetilde{W}\subset M;$
\footnote{ The  set $ \widetilde{W}$
is  defined   in Subsection \ref{subsection_cross} below.}

 \item[$\bullet$]
   for  every function $f$ as in the hypothesis of {\rm PROBLEM 2} with $Z=\C,$
 there exists a unique function
$\hat{f}\in\mathcal{O}(\widehat{\widetilde{W}}       \setminus \widehat{M},\C )$
such that
 $\hat{f}=f$ on $ (W\cap   \widetilde{W})\setminus M.$
  \end{itemize}
}

\smallskip

We refer the reader to \cite{jp3,jp5}  for complete versions of this theorem.

\smallskip

The  first result in the second direction (i.e.  ``boundary context'') is  contained in  the  work of  Malgrange--Zerner \cite{zen} in the
1960s.
Further  results in this direction were obtained by
Komatsu \cite{ko} and Dru\.{z}kowski \cite{dr}, but only for some special
cases.
Recently, Gonchar \cite{go1,go2} has proved  a  more
general result where the following case  of PROBLEM 1 has    been solved:  $D$ and $G$ are Jordan domains in $\C$,
  $ A$ (resp. $B$) is an open boundary  subset of $\partial D$ (resp. $\partial G$),   and   $ Z=\C.$
  Namely, we have

\medskip

\noindent {\bf Theorem  3
(Gonchar  \cite{go1,go2}).}
{\it
Let $X=Y=\C,$
let  $D\subset X,$ $G\subset Y$ be  Jordan domains and  $A $ (resp. $B$) a nonempty open set
of the boundary $\partial D$  (resp. $\partial G$).
Then,  for every  function  $f\in\mathcal{C}(W,\C)$  which  satisfies the hypotheses  of {\rm PROBLEM 1} with $Z=\C$,
 there exists a unique function
$\hat{f}\in\mathcal{C}(\widehat{W}\cup W,\C )      \cap \mathcal{O} (\widehat{W},\C)$
such that
 $\hat{f}=f$ on $ W.$
 Here
 \begin{equation*}
 \widehat{W}
 :=\left\lbrace (z,w)\in D\times G :\  \omega(z,A,D)+  \omega(w,B ,G)<1
\right\rbrace,
\end{equation*}
 where $ \omega(\cdot,A,D)$ and $ \omega(\cdot,B ,G)$ are the harmonic measures (see Subsection
 \ref{Subsection_approach_regions} below for this notion).
 }

 \smallskip

Theorem 3 may be rephrased as follows
 $\widehat{\widetilde{W}}= \widehat{W}$ (see also \cite{pn4}).
 It should be observed that before  Gonchar's works,  Airapetyan and Henkin published  a  version of
 the edge-of-the-wedge theorem for CR manifolds (see  \cite{ah1}  for a brief version
 and  \cite{ah2} for a complete proof).
   Gonchar's Theorem could be deduced from the latter result.
\section{New formulations} \label{sec_new_formulations}
Our purpose  is to develop a  theory which  unifies all  results  obtained so far.
First  we  develop  some  new notions  such  as  system of  approach  regions for an open set
in a  complex  manifold,
and  the corresponding plurisubharmonic measure. These  will provide the  framework
for an  exact  formulation of  PROBLEM 1 and 2,  and for  our solution.
 \subsection{Approach regions, local pluripolarity and plurisubharmonic measure}\label{Subsection_approach_regions}
\begin{defi}\label{defi_approach_region}
  Let $X$ be a complex manifold  and let $D\subset X$ be an open subset.
  A  {\rm system of approach regions} for $D$ is a collection
  $\mathcal{A}=\big(\mathcal{A}_{\alpha}(\zeta)\big)_{\zeta\in\overline{ D},\  \alpha\in I_{\zeta}}$
  ( $I_{\zeta}\neq\varnothing$ for all  $\zeta\in\partial D$) \footnote{ Note that this definition
  is  slightly different from Definition 2.1 in \cite{nv2}.}
  of open subsets of $D$
  with the following properties:
  \begin{itemize}
  \item[(i)] For all $\zeta\in D,$ the system $\big(\mathcal{A}_{\alpha}(\zeta)\big)_{ \alpha\in I_{\zeta}}$
  forms a basis of open neighborhoods of $\zeta$  (i.e., for any  open neighborhood $U$  of a point  $\zeta\in D,$
  there is $ \alpha\in I_{\zeta}$ such that  $\zeta\in \mathcal{A}_{\alpha}(\zeta)\subset U$).
  \item[(ii)] For  all $\zeta\in\partial D$ and
     $\alpha\in I_{\zeta},$   
  $\zeta\in \overline{\mathcal{A}_{\alpha}(\zeta)}.$
  \end{itemize}
 $\mathcal{A}_{\alpha}(\zeta)$ is often
  called an {\rm approach region} at $\zeta.$

 Moreover, $\mathcal{A}$  is
said to be {\rm canonical} if it satisfies (i) and the following  property (which is stronger than (ii)):
  \begin{itemize}
  \item[(ii')]
For every point  $\zeta\in \partial D,$
  there is  a basis of open neighborhoods  $(U_{\alpha})_{\alpha\in I_{\zeta}}$ of $\zeta$  in $X$ such that
   $ \mathcal{A}_{\alpha}(\zeta)=U_{\alpha}\cap D,$  $\alpha\in I_{\zeta}.$
   \end{itemize}
  \end{defi}

Various  systems of approach regions which one often encounters in Complex Analysis
 will be described in the  next  subsection.   Systems of approach regions  for $D$ are used to deal with the limit  at points in $\overline{D}$
 of  mappings defined on  some  open subsets of $D.$ Consequently, we deduce from
 Definition  \ref{defi_approach_region} that
 the subfamily $\big(\mathcal{A}_{\alpha}(\zeta)\big)_{\zeta\in D,\ \alpha\in I_{\zeta}}$
 is, in  a certain sense, independent of the choice of a system of approach regions  $\mathcal{A}.$
 In addition, any two canonical systems of approach regions are, in some sense, equivalent.
 These observations lead us to use, throughout the paper,  the following {\bf convention}:

{\it We fix, for every open set $D\subset X,$  a {\bf canonical system of approach regions.}
 When we want to define a  system of approach regions $\mathcal{A}$ for an open set $D\subset X,$
 we only need to specify the subfamily   $\big(\mathcal{A}_{\alpha}(\zeta)\big)_{\zeta\in\partial
  D,\ \alpha\in I_{\zeta}}.$}

  In what follows we fix an open subset $D\subset X$ and a  system of approach regions   $\mathcal{A}=\big(\mathcal{A}_{\alpha}(\zeta)\big)_{
  \zeta\in\overline{ D},\
   \alpha\in I_{\zeta}}$ for $D.$

 For every function $u:\ D\longrightarrow [-\infty,\infty),$ let
\begin{equation*}
 (\Alimsup u)(z):=
 \sup\limits_{\alpha\in I_{z}}\limsup\limits_{w\in \mathcal{A}_{\alpha}(z),\ w\to z}u(w), \qquad  z\in\overline{D}.
\end{equation*}
By  Definition \ref{defi_approach_region} (i),
 $(\Alimsup u)|_D$ coincides with the usual
{\it upper semicontinuous regularization}  of $u.$

For a set  $A\subset \overline{D}$ put
\begin{equation*}
h_{A,D}:=\sup\left\lbrace u\ :\  u\in\mathcal{PSH}(D),\ u\leq 1\ \text{on}\ D,\
   \Alimsup u\leq 0\ \text{on}\ A    \right\rbrace,
\end{equation*}
where $\mathcal{PSH}(D)$ denotes the cone  of all functions  plurisubharmonic
on $D.$

 $A\subset D$ is
said to be  {\it thin} in $D$ if  for    every point $a\in D$  there  is a connected
 neighborhood $U=U_a\subset  D$ and  a holomorphic  function $f$ on $U,$ not identically  zero such that  $U\cap A\subset f^{-1}(0).$
   $A\subset D$ is said to be {\it pluripolar} in $D$ if there is  $u\in
\mathcal{PSH}(D)$ such that $u$ is not identically $-\infty$ on every
connected component of $D$ and $A\subset \left\lbrace z\in D:\ u(z)=-\infty\right\rbrace.$
$A\subset D$ is said to be {\it locally  pluripolar} in $D$ if  for any
$z\in A,$ there is an open  neighborhood $V\subset D$ of $z$ such that $A\cap V$ is
pluripolar in $V.$ $A\subset D$ is said  to be {\it nonpluripolar} (resp. {\it non locally  pluripolar}) if it is not
pluripolar (resp. not locally pluripolar).
According to a classical result of Josefson and Bedford  (see \cite{jo},
\cite{be}), if $D$ is a Riemann--Stein domain 
then   $A\subset D$ is   locally  pluripolar if and only if it is
pluripolar.


\begin{defi}\label{defi_relative_extremal}
 For  $A\subset\overline{D},$ the
  {\rm relative extremal function of $A$  relative to $D$} is
 the function $ \omega(\cdot,A,D)$
 defined by
 \begin{equation*}
 \omega(z,A,D)= \omega_{\mathcal{A}}(z,A,D):= (\Alimsup  h_{A,D})(z),\qquad  z\in \overline{D}.\quad \footnotemark
\end{equation*}
\end{defi}
\footnotetext{Observe that this function depends on the system of
approach regions.}
 Note that when $A\subset D,$  Definition \ref{defi_relative_extremal}
   coincides with the classical
definition of Siciak's  relative extremal function. When  $D$ is  a  complex manifold of dimension $1$ and  $\mathcal{A}$  is the
canonical  system,  the function $ \omega(\cdot,A,D)$ is often called  {\it the harmonic  measure of $A$  relative to $D$}  (see Theorem 3 above).

 Next, we say that a  set  $A\subset \overline{D}$ is {\it
locally pluriregular at a point $a\in \overline{A}$}   if $\omega(a,A\cap U,D\cap U)=0$
  for  all open neighborhoods $U$ of $a.$
Moreover, $ A$ is said to be {\it locally pluriregular } if it is locally
pluriregular at all points $a\in A.$
It should be noted from  Definition  \ref{defi_approach_region} that  if $a\in \overline{A}\cap D$ then  the property of local pluriregularity
of $A$ at $a$  does not depend  on any particular choices of  a system of approach regions $\mathcal{A},$
while the situation is different when  $a\in \overline{A}\cap \partial D:$ the  property    does depend on  $\mathcal{A}.$

 We denote by $A^{\ast}$ the following set
 \begin{equation*}
  (A\cap\partial D)\bigcup\left\lbrace   a\in \overline{A}\cap D:\ A\ \text{is locally pluriregular at}\ a \right\rbrace.
  \end{equation*}
If $A\subset D$ is non locally pluripolar, then  a classical result of Bedford and Taylor
(see \cite{be,bt})
says that    $A^{\ast}$ is locally
pluriregular    
 and  $A\setminus A^{\ast}$ is locally  pluripolar.
Moreover, $A^{\ast}$ is locally  of type $\mathcal{G}_{\delta},$  that  is,
for every $a\in A^{\ast}$ there is an open  neighborhood $U\subset D$ of $a$
such that $A^{\ast}\cap U$ is  a countable intersection of open sets.

Now we are  in the position to formulate the following  version
of the  plurisubharmonic measure.
\begin{defi}\label{defi_pluri_measure}
For a set $A\subset \overline{D},$ let
$\widetilde{A} =\widetilde{A}(\mathcal{A}):=\bigcup\limits_{P\in \mathcal{E}(A)} P,$ where
\begin{equation*}
\mathcal{E}(A)=\mathcal{E}(A,\mathcal{A}):=\left\lbrace P\subset \overline{D}:\  P\ \text{is  locally pluriregular,}\
\overline{P}\subset  A^{\ast}      \right\rbrace,
\end{equation*}
 The {\rm  plurisubharmonic measure of $A$  relative to $D$} is
 the function $\widetilde{\omega}(\cdot,A,D)$
 defined by
\begin{equation*}
\widetilde{\omega}(z,A,D):=  \omega(z,\widetilde{A},D),\qquad  z\in D.
\end{equation*}
\end{defi}

It is worthy to remark that  $\widetilde{\omega}(\cdot,A,D)\in\mathcal{PSH}(D)$ and
$0\leq \widetilde{\omega}(z,A,D)\leq 1,\ z\in  D.$  Moreover,
\begin{equation}\label{eq_defi_pluri_measure}
\Big(\Alimsup \widetilde{\omega}(\cdot,A,D)\Big)(z)=0,\qquad z\in \widetilde{A}.
\end{equation}
 An example  in  \cite{ah} shows  that,     in general,
  $\omega(\cdot,A,D)\not=\widetilde{\omega}(\cdot,A,D)$ on $D.$
Section  \ref{Section_PROBLEM_1_second_case}
and  \ref{section_application} below are devoted  to the study of  $ \widetilde{\omega}(\cdot,A,D)$
in some important cases.  As we will see later,  in most applications  one can obtain    good and  simple
characterizations  of      $ \widetilde{\omega}(\cdot,A,D)$  (see Theorem 5, 6, 7,  9 and Corollary 2, 3  below).

Now  we compare the plurisubharmonic measure $ \widetilde{\omega}(\cdot,A,D)$ with
Siciak's relative  extremal function $ \omega(\cdot,A,D).$   We only  consider two
important special cases:  $A\subset D$ and  $A\subset\partial D.$
For the moment, we only focus  on the case  where  $A\subset D.$
The latter one  will be discussed  in
Section  \ref{Section_PROBLEM_1_second_case}  and \ref{section_application} below.

  If $A$ is an open subset of an arbitrary complex manifold  $D$, then it can be shown that
\begin{equation*}
\widetilde{\omega}(z,A,D)=  \omega(z,A,D),\qquad  z\in D.
\end{equation*}
If $A$  is  a (not necessarily open) subset of an arbitrary complex manifold $D,$
then we have,
by Proposition 7.1 in  \cite{nv2},
  \begin{equation*}
\widetilde{\omega}(z,A,D)=  \omega(z,A^{\ast},D),\qquad  z\in D.
\end{equation*}
On the other hand,  if, morever,  $D$  is  a  bounded  open subset of  $\C^n$   then we  have (see, for  example,
  Lemma 3.5.3 in \cite{jp1}) $
\omega(z,A,D)=  \omega(z,A^{\ast},D),$ $ z\in D.$
Consequently, under  the last assumption,
\begin{equation*}
\widetilde{\omega}(z,A,D)=  \omega(z,A,D),\qquad  z\in D.
\end{equation*}
 Our discussion  shows that at least  in the case  where  $A\subset D$, the notion of  the plurisubharmonic measure
is  a  good candidate  for  generalizing  Siciak's relative  extremal function  to the manifold context in the theory
of separate holomorphy.

 For a good background of the pluripotential
theory, see the books  \cite{jp1} or  \cite{kl}.
\subsection{Examples of systems of  approach  regions} \label{Subsection_Examples}
  There  are  many  systems    of  approach  regions which are very useful
in Complex  Analysis. In this  subsection  we present some of them.

\noindent {\bf 1. Canonical  system   of  approach regions.}
It has been given by   Definition \ref{defi_approach_region} (i)--(ii').
 This is  the most natural one.

\noindent {\bf 2. System of  angular (or  Stolz)    approach regions for the open unit  disc.}
Let $E$  be the  open unit  disc of  $\C.$ Put
\begin{equation*}
\mathcal{A}_{\alpha}(\zeta):=
 \left\lbrace          t\in E:\ \left\vert
 \arg\left(\frac{\zeta-t}{\zeta}\right)
 \right\vert<\alpha\right\rbrace   ,\qquad  \zeta\in\partial E,\  0<\alpha<\frac{\pi}{2},
\end{equation*}
where  $\arg:\ \C\longrightarrow (-\pi,\pi]$ is  as usual the argument function.
$ \mathcal{A}=\left(\mathcal{A}_{\alpha}(\zeta)\right)_{\zeta\in\partial E,\
 0<\alpha<\frac{\pi}{2}}$ is   referred to as   {\it   the system of  angular (or  Stolz)
    approach regions for $E.$}
In this  context  $\Alim$ is also called {\it angular limit}.

\noindent {\bf 3. System of   angular    approach regions for   certain ``good" open subsets of Riemann surfaces.}
Now we  generalize the previous construction (for the open unit disc) to  a  global  situation. More  precisely, we
will use as  the  local  model the  system    of   angular    approach regions for  $E.$
Let $X$ be  a  complex  manifold of  dimension $1$ (in other words, $X$ is  a Riemann surface),
and  let $D\subset X$ be an  open  set. Then $D$   is   said  to be  {\it
good at a  point $\zeta\in\partial D$}  \footnote{ In the  work \cite{pn2}
 we  use  the more  appealing  word  {\it Jordan-curve-like} for this  notion.}
 if there is a  Jordan domain
   $U\subset X$  such that   $\zeta\in U$ and $U\cap
\partial D$ is the  interior of  a  Jordan curve.

Suppose that   $D$   is  good  at $\zeta.$ This point is   said  to be
 {\it of type 1} if there  is a neighborhood $V$ of
$\zeta$ such that $V_0=V\cap D$ is  a Jordan domain.  Otherwise,
$\zeta$ is   said  to be {\it   of type 2}. We  see  easily  that if
$\zeta$ is of  type 2, then  there  are an open  neighborhood    $V$ of $\zeta$ and  two disjoint  Jordan domains  $V_1,$ $V_2$
 such that  $V\cap
D=V_1\cup V_2.$  Moreover,    $D$  is   said  to be   {\it
good  on  a  subset  $A$}  of  $\partial D$ if  $D$ is
good at all points  of   $A.$

Here  is  a simple  example   which may clarify the  above definitions.
Let $G$ be  the  open  square in $\C$
with vertices $1+i,$ $-1+i,$ $-1-i,$ and  $1-i.$
Define the domain
\begin{equation*}
D:=G\setminus \left [-\frac{1}{2},\frac{1}{2}\right].
\end{equation*}
Then   $D$ is  good on $\partial G\cup  \left
(-\frac{1}{2},\frac{1}{2}\right).$ All points  of $\partial G$ are of  type
1 and all  points of $\left(-\frac{1}{2},\frac{1}{2}\right)$ are of type
2.

Suppose  now  that   $D$  is  good
on a  nonempty subset  $A$ of $\partial D.$ 
 We define {\it the system    of angular approach regions   supported on  $A$:}
   $\mathcal{A}=\big(\mathcal{A}_{\alpha}(\zeta)\big)_{\zeta\in\overline{ D},\  \alpha\in I_{\zeta}}$
as  follows:
  \begin{itemize}
  \item[$\bullet$] If $\zeta\in  \overline{ D}\setminus A,$  then $\big(\mathcal{A}_{\alpha}(\zeta)\big)_{ \alpha\in I_{\zeta}}$
coincide with the canonical approach  regions.
 \item[$\bullet$] If $\zeta\in  A,$ then  by using  a conformal  mapping $\Phi$ from $V_0$  (resp. $V_1$ and $V_2$)
 onto $E$
 when $\zeta$ is  of  type 1  (resp. 2), we can ``transfer'' the  angular  approach regions
at the  point $\Phi(\zeta)\in \partial E:$ $\left(\mathcal{A}_{\alpha}(\Phi(\zeta))\right)_{
 0<\alpha<\frac{\pi}{2}}$ to those at the point $\zeta\in  \partial D$  (see \cite{pn2} for more  detailed explanations).
\end{itemize}
Making use  of conformal  mappings in  a  local  way, we  can   transfer, in  the  same  way, many
notions which  exist   on  $E$  (resp.  $\partial E$)  to those  on  $D$  (resp.  $\partial D$).

\noindent {\bf 4. System     of conical approach regions.}

Let $D\subset\C^n$ be a  domain  and $A\subset\partial D.$ Suppose in  addition that
for every  point $\zeta\in A$ there exists the (real) tangent space $T_{\zeta}$ to $\partial D$ at $\zeta.$
 We  define {\it the system   of   conical approach regions supported on $A$:}
   $\mathcal{A}=\big(\mathcal{A}_{\alpha}(\zeta)\big)_{\zeta\in\overline{ D},\  \alpha\in I_{\zeta}}$
as  follows:
  \begin{itemize}
  \item[$\bullet$] If $\zeta\in  \overline{ D}\setminus A,$  then $\big(\mathcal{A}_{\alpha}(\zeta)\big)_{ \alpha\in I_{\zeta}}$
coincide with the  canonical    approach regions.
 \item[$\bullet$] If $\zeta\in  A,$ then
 \begin{equation*}
 \mathcal{A}_{\alpha}(\zeta):=\left\lbrace  z\in D:\  \vert  z-\zeta\vert <\alpha\cdot \dist(z,T_{\zeta})  \right\rbrace,
 \end{equation*}
where $I_{\zeta}:=(1,\infty)$  and $ \dist(z,T_{\zeta})$ denotes the Euclidean    distance  from the  point $z$ to $T_{\zeta}.$
\end{itemize}

We can also   generalize the previous construction to  a  global situation:

{\it  $X$ is an  arbitrary complex manifold, $D\subset  X$ is an open set and  $A\subset \partial D$ is  a  subset
with  the property that
at every point $\zeta\in A$
 there exists the  (real)  tangent  space  $T_{\zeta}$ to $\partial D.$}

We can  also    formulate  the  notion of points of  type 1  or  2 in this  general  context  in the same way
as we have  already done
in Paragraph 3  above of this  subsection.

\subsection{Cross and  separate holomorphicity and $\mathcal{A}$-limit.} \label{subsection_cross}
Let $X,\ Y$  be two complex manifolds,
  let $D\subset X,$ $ G\subset Y$ be two nonempty open sets, let
  $A\subset \overline{D}$   and  $B\subset \overline{G}.$
  Moreover, $D$  (resp.  $G$) is equipped with a
  system of approach regions
  $\mathcal{A}(D)=\big(\mathcal{A}_{\alpha}(\zeta)\big)_{\zeta\in\overline{D},\  \alpha\in I_{\zeta}}$
  (resp.  $\mathcal{A}(G)=\big(\mathcal{A}_{\alpha}(\eta)\big)_{\eta\in\overline{G},\  \alpha\in I_{\eta}}$).\;\footnote{ In fact we should have written $I_\zeta(D)$,
  resp. $I_\eta(G)$; but we skip $D$ and $G$ here to make the notions
  as simple as possible.}
 We define
a {\it $2$-fold cross} $W,$  its {\it  interior} $W^{\text{o}}$ and
its {\it  regular part} $\widetilde{W}$ (with respect to
$\mathcal{A}(D)$ and $\mathcal{A}(G)$)
as
\begin{eqnarray*}
W &=&\X(A,B; D,G)
:=\big((D\cup A)\times B\big)\bigcup\big (A\times(B\cup G)\big),\\
W^{\text{o}} &=&\X^{\text{o}}(A,B; D,G)
:= (A\times  G)\cup (D\times B),\\
\widetilde{W} &=&\widetilde{\X}(A,B;D,G) := \big((D\cup
\widetilde{A})\times \widetilde{B}\big)\bigcup
\big(\widetilde{A}\times(G\cup \widetilde{B})\big),
\end{eqnarray*}
where  $\widetilde{A}$ and $ \widetilde{B}$ are calculated  using
Definition \ref{defi_pluri_measure}.
Moreover, put
\begin{eqnarray*}
\omega(z,w)&:=&\omega(z,A,D)+\omega(w,B,G),\qquad
(z,w)\in D\times G,\\
\widetilde{\omega}(z,w)&:=&\widetilde{\omega}(z,A,D)+\widetilde{\omega}(w,B,G),\qquad
(z,w)\in D\times G.
\end{eqnarray*}

For a $2$-fold cross $W :=\X(A,B; D,G)$ let
\begin{eqnarray*}
\widehat{W}&:=&\widehat{\X}(A,B;D,G) =\left\lbrace (z,w)\in D\times
G:\ \omega(z,w)  <1
\right\rbrace,\\
 \widehat{\widetilde{W}} &:=&\widehat{\X}(\widetilde{A},\widetilde{B};D,G)
 =\left\lbrace (z,w)\in D\times G :\  \widetilde{\omega}(z,w)<1
\right\rbrace.
\end{eqnarray*}

Let $Z$ be a complex analytic space and $M\subset W$ a subset which is  relatively  closed in fibers over $A$ and $B.$
 We say that a mapping
$f:W^{\text{o}}\setminus M \longrightarrow Z$ is {\it separately holomorphic}
  and write $f\in\mathcal{O}_s(W^{\text{o}}\setminus M ,Z),$   if,
 for any $a\in A $ (resp.  $b\in B$)
 the mapping $f(a,\cdot)|_{G\setminus M_a}$  (resp.  $f(\cdot,b)|_{D\setminus M^b}$)  is holomorphic.

 We say that a mapping $f:\  W\setminus M \longrightarrow Z$
   is  {\it separately continuous}
and write
 $f\in \mathcal{C}_s\Big( W\setminus M,Z  \Big)$
 if,
 for any $a\in A$ (resp.  $b\in B$)
 the mapping $f(a,\cdot)|_{(G\cup B)\setminus M_a   }$  (resp.  $f(\cdot,b)|_{(D\cup A)\setminus  M^b}$)  is continuous.

Let  $\Omega$ be  an  open subset of $D\times  G.$ A  point $(\zeta,\eta)\in\overline{D}\times\overline{G}$ is  said
to be an {\it end-point} of $\Omega$  with  respect to
 $\mathcal{A}=\mathcal{A}(D)\times \mathcal{A}(G)$ if for any $(\alpha,\beta)\in I_{\zeta}\times I_{\eta}$ there
 exist  open neighborhoods $U$ of $\zeta$
 in $X$ and $V$ of $\eta$ in $Y$   such that
  \begin{equation*}
  \Big(U\cap \mathcal{A}_{\alpha}(\zeta)\Big)   \times \Big(V\cap \mathcal{A}_{\beta}(\eta)
  \Big)\subset \Omega.
  \end{equation*}
 The  set  of  all end-points  of $\Omega$  is  denoted  by  $\End(\Omega).$

 It follows  from  (\ref{eq_defi_pluri_measure}) that if $\widetilde{A},\widetilde{B}\not=\varnothing$,
  then $  \widetilde{W}\subset \End(
 \widehat{\widetilde{W}}).$

Let $S$ be  a relatively closed  subset of $\widehat{\widetilde{W}}$  and  let $
(\zeta,\eta)\in\End(\widehat{\widetilde{W}}\setminus  S).$    Then a mapping $f:\ \widehat{\widetilde{W}}\setminus
S\longrightarrow Z$ is said to {\it admit the $\mathcal{A}$-limit $\lambda$ at $ (\zeta,\eta),$}  and one writes
\begin{equation*}
(\Alim f)(\zeta,\eta)=\lambda,\qquad\footnote{ Note that here $\mathcal{A}=\mathcal{A}(D)\times \mathcal{A}(G).$}
\end{equation*}
 if, for all $\alpha\in I_{\zeta},\  \beta\in I_{\eta},$
\begin{equation*}
\lim\limits_{ \widehat{\widetilde{W}}\setminus S\ni (z,w)\to (\zeta,\eta),\ z\in \mathcal{A}_{\alpha}(\zeta),\  w\in
\mathcal{A}_{\beta}(\eta)}f(z,w)=\lambda.
\end{equation*}

We conclude this introduction with a notion we need in the sequel.  Let   $\mathcal{M}$ be  a topological  space.
 A
mapping  $f:\ \mathcal{M}\longrightarrow  Z$ is said to be {\it
bounded} if there exists an open neighborhood $U$ of
$f(\mathcal{M})$ in $Z$ and a holomorphic embedding $\phi$ of $ U $
into the unit  polydisc of $ \C^k$  such that $\phi(U)$ is an analytic set
in this polydisc.
      $f$ is said to be {\it locally bounded along} $\mathcal{N}\subset \mathcal{M}$ if
for every point $z\in \mathcal{N},$ there is an open neighborhood
$U$ of $z$ (in $\mathcal{M}$) such that  $f|_{U}:\ U \longrightarrow
Z$ is bounded.   $f$ is said to be {\it locally bounded} if  it is
so for  $\mathcal{N}= \mathcal{M}.$ It is clear that,  if $Z=\C$,
then the above notions of boundedness coincide with the usual ones.
%
%
%
\subsection{Hartogs extension property.}\label{HEP}
The following example (see Shiffman  \cite{sh2}) shows that an additional
hypothesis on the ``target space" $Z$ is necessary in order that   PROBLEM 1 and 2  make sense. Consider the mapping $f:\
\C^2\longrightarrow \P^1$ given by
\begin{equation*}
  f(z,w):=
\begin{cases}
[(z+w)^2: (z-w)^2],
  & (z,w)\not=(0,0) ,\\
 [1:1], &  (z,w)=(0,0).
\end{cases}
\end{equation*}
Then $f\in\mathcal{O}_s\Big(\X^{\text{o}}(\C,\C;\C,\C),\P^1\Big),$ but $f$ is not
continuous at $(0,0).$

 We recall here  the following notion
(see,  for  example, Shiffman \cite{sh1}).
 Let $p\geq 2$ be an integer. For $0<r<1,$  the {\it Hartogs
 figure} in dimension $p,$ denoted by $H_p(r),$ is given by
 \begin{equation*}
H_p(r):=\left\lbrace (z^{'},z_p)\in E^p: \  \Vert z^{'}\Vert<r \ \ \text{or}\ \ \vert z_p\vert >1-r \right\rbrace,
 \end{equation*}
where $E$ is the open  unit disc of $\C$ and $z^{'}=(z_1,\ldots,z_{p-1}),$
$\Vert z^{'}\Vert:=\max\limits_{1\leq j\leq p-1} \vert z_j\vert.$

\begin{defi}\label{defi_HEP}
A complex analytic space $Z$ is said to {\rm possess the Hartogs extension property
  in dimension $p$} if  every  mapping $f\in\mathcal{O}(H_p(r) , Z)$ extends to a  mapping
  $\hat{f}\in\mathcal{O}(E^p , Z).$  Moreover, $Z$ is said to {\rm possess the Hartogs extension property}
 if it possesses this  property in all dimensions $p\geq 2.$
\end{defi}
It is a classical result of Ivashkovich (see \cite{iv1}) that if $Z$ possesses
the Hartogs extension property  in dimension 2, then it possesses this property in all dimensions $p\geq 2.$
 Some typical  examples of complex analytic  spaces  possessing the Hartogs extension property are  the
complex Lie groups (see \cite{asy}), the taut spaces (see \cite{wu}), the Hermitian manifold with
negative holomorphic
sectional  curvature (see \cite{sh1}), the holomorphically convex K\"{a}hler  manifold without rational curves
(see \cite{iv1}).

Here we  mention  an important characterization.

\medskip

\noindent {\bf Theorem 4  (Shiffman \cite{sh1}).}
  {\it
A complex analytic space $Z$   possesses the Hartogs extension property
    if and only if  for every domain $D$ of any Stein manifold $\mathcal{M},$ every mapping
    $f\in\mathcal{O}(D, Z)$ extends to a  mapping $\hat{f}\in
    \mathcal{O}(\widehat{D},Z) ,$   where $\widehat{D}$ is the envelope of holomorphy\footnote{ For the notion of the envelope of holomorphy, see, for example, \cite{jp1}.}
of $D.$
}

\medskip

In the light of Definition  \ref{defi_HEP} and   Shiffman's Theorem,
 the natural ``target spaces" $Z$ for obtaining
satisfactory answers to  PROBLEM 1 are the complex analytic spaces  which possess the Hartogs extension property.

\section{A new approach: Poletsky Theory of discs and Rosay Theorem}  \label{section_new_approach}

 Poletsky Theory of discs was invented by Poletsky (see \cite{po1,po2}) at the end of the 1980s.
A new approach  to the theory of
separate holomorphy based on   Poletsky theory of discs  was developed     in our work \cite{nv1}.
Let us recall some elements of this theory.

 Let $E$ denote as usual the open unit disc in $\C.$ %
For a complex
manifold $\mathcal{M},$ let $\mathcal{O}(\overline{E},\mathcal{M})$ denote
the set of all holomorphic mappings $\phi:\ E\longrightarrow \mathcal{M}$ which
extend holomorphically  to   a neighborhood of  $\overline{E}.$
Such a mapping $\phi$ is called a {\it holomorphic disc} on $\mathcal{M}.$ Moreover, for
a subset $A$ of $\mathcal{M},$ let
\begin{equation*}
 1_{  A,\mathcal{M}}(z):=
\begin{cases}
1,
  &z\in   A,\\
 0, & z\in \mathcal{M}\setminus A.
\end{cases}
\end{equation*}

In 2003  Rosay proved the following remarkable result.

\renewcommand{\thethmspec}{Rosay Theorem   (\cite{ro})}
  \begin{thmspec}\label{Rosaythm}
Let $u$ be an upper semicontinuous function on a complex manifold
$\mathcal{M}.$ Then the Poisson functional of $u$  defined by
\begin{equation*}
\mathcal{P}[u](z):=\inf\left\lbrace\frac{1}{2\pi}\int\limits_{0}^{2\pi} u(\phi(e^{i\theta}))d\theta:  \
\phi\in   \mathcal{O}(\overline{E},\mathcal{M}), \ \phi(0)=z
\right\rbrace,
\end{equation*}
is plurisubharmonic on $\mathcal{M}.$
\end{thmspec}

Rosay Theorem may be viewed as an important development in Poletsky
theory of   discs. Observe that special cases of this  theorem
 have been considered by Poletsky (see \cite{po1,po2}),
L\'arusson--Sigurdsson (see \cite{ls}) and Edigarian (see \cite{ed}).

 The next result  describes the situation  in dimension $1.$
 \renewcommand{\thethmspec}{Lemma 1 (\cite[Lemma 3.3]{nv1})}
  \begin{thmspec}\label{lem_circle}
Let $T$ be an open subset of $\overline{E}.$   Then
\begin{equation*}
\omega(0,T\cap E,E)\leq \frac{1}{2\pi}\int\limits_{0}^{2\pi}1_{  \partial E\setminus T,T}
(e^{i\theta})d\theta.
\end{equation*}
\end{thmspec}
The last result, which is an important consequence of Rosay Theorem,  gives the connection between
  the Poisson functional       and the plurisubharmonic measure.
\renewcommand{\thethmspec}{Lemma 2 (\cite[Proposition 3.4]{nv1})}
  \begin{thmspec}\label{prop_Rosay}
Let $\mathcal{M}$ be a complex manifold and   $A$   a nonempty open
subset of $\mathcal{M}.$    Then
   $ \omega(z,A,\mathcal{M}) = \mathcal{P}[1_{\mathcal{M}\setminus A,\mathcal{M}}](z),$ $z\in\mathcal{M}.$
\end{thmspec}

%
%
%
\section{PROBLEM 1 for the case $A\subset D,$ $B\subset G$}
\label{Section_PROBLEM_1_first_case}
We will give  the first application of   the previous section.
Observe that  under the hypothesis $A\subset D,$ $B\subset G$  and  the notation
of Subsection  \ref{subsection_cross}, we have  $W=W^{\text{o}}$
and $W\cap \widetilde{W}\subset  W\cap\widehat{\widetilde{W}}.$
Since  $ \widetilde{W}\subset D\times G,$  the notion $\Alim$  at a point of $ \widetilde{W}$
 coincides with the ordinary notion of a limit, that is, $\mathcal{A}$  can be  taken as the canonical  system. Moreover, it can be  shown that
$W\setminus  \widetilde{W}$ is a locally pluripolar subset of $D\times G.$
Therefore,  from the viewpoint of the pluripotential theory, $W\cap  \widetilde{W}$
is ``almost"  equal to  $W.$
Now  we are able to state   the following  generalization of Theorem 1.
\renewcommand{\thethmspec}{Theorem 5   (\cite[Theorem A]{nv1})}
  \begin{thmspec}
  Let $X,\ Y$ be   arbitrary complex manifolds, let  $D\subset X$ and $G\subset Y$ be     open sets
  and  $ A\subset
D,$  $B\subset G$  non locally pluripolar   subsets.
 Let $Z$ be a complex  analytic space possessing the
  Hartogs
  extension property.
Then for every mapping $f\in\mathcal{O}_s(W^{\text{o}},Z),$ there is a unique mapping
$\hat{f}\in\mathcal{O}(\widehat{\widetilde{W}},Z)$
such that $\hat{f}=f$ on $W\cap \widehat{\widetilde{W}}.$
\end{thmspec}
A remark is in order. Theorem  5 removes all the assumptions
 of pseudoconvexity of the ``source spaces" $X,$ $Y$
stated in  Theorem 1.
 Namely, now $X$ and $Y$ can be  arbitrary complex manifolds.
  The sketchy proof given below explain  our approach:  how  Poletsky theory of discs
  and Rosay Theorem
may apply to the theory of separate holomorphy.
It is divided into four steps.  In Step 3 and 4 below  we use some ideas in our previous
joint-work with Pflug  \cite{pn1}.

\smallskip

\noindent {\bf Step 1:} {\it  The case
where   $D$ is an arbitrary complex manifold, $A$ is  an open subset of $D,$ and    $G$ is a bounded open subset of $\C^n.$}

\smallskip

\noindent{ \it Sketchy proof of Step 1.} We define $\hat{f}$ as follows: Let $\mathcal{W}$ be  the set of all pairs
$(z,w)\in D\times G$  with the property that there are a holomorphic disc
$\phi\in\mathcal{O}(\overline{E},D)$ and $t\in E$ such that $\phi(t)=z$
and  $(t,w)\in\widehat{\widetilde{\X}}\left(\phi^{-1}(A)\cap E,B;E,G\right).$  In virtue of
 Theorem 1 and the observation made at the beginning of the section,
 let $\hat{f}_{\phi}$ be the unique mapping  in
$ \mathcal{O}\left(\widehat{\widetilde{\X}}(\phi^{-1}(A)\cap E,B;E,G),Z\right)$ such that
\begin{equation} \label{eq1_CaseI}
  \hat{f}_{\phi}(t,w)=f(\phi(t),w),\qquad (t,w)\in \X\left(\phi^{-1}(A)\cap E,B;E,G\right)
  \cap \widetilde{\X}\left(\phi^{-1}(A)\cap E,B;E,G\right).
\end{equation}
   Then we may define the desired extension mapping $\hat{f}$ as  follows
\begin{equation}\label{eq2_CaseI}
 \hat{f}(z,w):=\hat{f}_{\phi}(t,w) .
\end{equation}
Using the uniqueness of Theorem 1, we can prove that $\hat{f}$ is well-defined on $\mathcal{W}.$
Using Lemma 1 and 2, one can  show that
\begin{equation*}
 \mathcal{W}= \widehat{\widetilde{W}}  .
\end{equation*}
 Moreover, it follows from the above construction  that for every
fixed $z\in D,$  the restricted mapping $\hat{f}(z,\cdot)$ is holomorphic on the open set
$\left\lbrace w\in G:\ (z,w)\in \widehat{\widetilde{W}} \right\rbrace.$
However,  it  is  quite difficult to see that   $\hat{f}$ is holomorphic in both variables $(z,w).$
A complete proof of this fact is  given in Theorem 4.1 in \cite{nv1}.  Now  we only  explain   briefly  why $\hat{f}$ is holomorphic in a neighborhood of an arbitrarily fixed point $(z_0,w_0)\in \widehat{\widetilde{W}} .$ For this purpose  we  ``add" one complex dimension more to a suitable  neighborhood of $(z_0,w_0),$ and this makes
our initial $2$-fold cross $W$ a  $3$-fold one. Finally, we try to apply  the version of Theorem 1  for $3$-fold cross in order
to finish the proof.
\hfill $\square$

\smallskip

\noindent {\bf Step 2:} {\it  The case
where
   $D,\ G$ are arbitrary complex manifolds,   but $A\subset D,$ $B\subset G$ are  open subsets.}

   \smallskip

\noindent{ \it Sketchy proof of Step 2.}
 It follows from the discussion made at the end of Subsection \ref{Subsection_approach_regions}
  that
 under the hypothesis of Step 2,  $\widehat{\widetilde{W}}=\widehat{W}$
 and $W=W\cap\widetilde{W}\subset \widehat{W}.$

We will determine  the value of $\hat{f}$ at   an arbitrary fixed point $(z_0,w_0)\in\widehat{W}.$
To this end fix any $\epsilon>0$ such that
\begin{equation}\label{eq_prop2}
2\epsilon<1- \omega(z_0,A,D) - \omega(w_0,B,G).
\end{equation}
 Applying Rosay Theorem   and Lemma 2, there is a holomorphic disc   $\phi\in
\mathcal{O}(\overline{E} , D)$  (resp.  $\psi\in\mathcal{O}(
\overline{E} , G))$  such that $\phi(0)=z_0$  (resp. $\psi(0)=w_0$) and
\begin{equation*}
\frac{1}{2\pi} \int\limits_{0}^{2\pi} 1_{D\setminus
A}(\phi(e^{i\theta}))d\theta< \omega(z_0,A,D)+\epsilon,\qquad
\frac{1}{2\pi} \int\limits_{0}^{2\pi} 1_{G\setminus
B}(\psi(e^{i\theta}))d\theta< \omega(w_0,B,G)+\epsilon.
\end{equation*}
Using this and estimate (\ref{eq_prop2}) and Lemma 1,  we see that
  \begin{equation*}
  (0,0)\in\widehat{\X}\left(\phi^{-1}(A)\cap E,\psi^{-1}(B)\cap E;E,E\right)
 .
 \end{equation*}
  Moreover, since $f\in\mathcal{O}_s(W^{\text{o}},Z),$ the mapping $h$ given by
 \begin{equation*}
  h(t,\tau):= f(\phi(t),\psi(\tau)),\qquad (t,\tau)\in
  \X\left(\phi^{-1}(A)\cap E,\psi^{-1}(B)\cap E;E,E\right),
 \end{equation*}
 belongs to $\mathcal{O}_s\Big(\X\left(\phi^{-1}(A)\cap E,\psi^{-1}(B)\cap E;E,E\right),Z\Big).$
  By Theorem  1, let $\hat{h}\in \mathcal{O}\left(\widehat{\X}\left(\phi^{-1}(A)\cap E,\psi^{-1}(B)\cap E;E,E
 \right),Z\right)$  be the unique  mapping  such that
 \begin{equation*}
 \hat{h}(t,\tau)=h(t,\tau)= f(\phi(t),\psi(\tau)),\qquad
(t,\tau)\in \X\left(\phi^{-1}(A)\cap E,\psi^{-1}(B)\cap E;E,E\right).
\end{equation*}
Then we can define
\begin{equation*}
\hat{f}(z_0,w_0)=\hat{h}(0,0),\qquad  (z_0,w_0)\in\widehat{W}.
\end{equation*}
We leave to the interested reader the verification that  $\hat{f}$ is well-defined on $ \widehat{W}.$
Now we explain why  $\hat{f}\in\mathcal{O}(  \widehat{W},Z).$

If we fix  $\phi$ and let $\psi$ be free (or conversely, fix   $\psi$ and let $\phi$ be free)
in the above construction,  then this procedure is   very similar to the one
carried out in (\ref{eq1_CaseI})--(\ref{eq2_CaseI}). Consequently, we may apply the result of Step  1
twice in order to conclude that for all $(z_0,w_0)\in  \widehat{W},$
 $\hat{f}(z_0,\cdot)$ is holomorphic in  $\{ w\in G:\ (z_0,w)\in \widehat{W} \}$
(resp.      $\hat{f}(\cdot,w_0)$ is holomorphic in  $\{ z\in D:\ (z,w_0)\in \widehat{W} \}$).
Applying the classical Hartogs extension theorem, it follows that  $\hat{f}\in\mathcal{O}(\widehat{W},Z).$
\hfill $\square$

\smallskip

 To continue  the proof  we need to introduce some more notation.

Suppose without loss of generality that $D$ and $G$ are domains and let $m$ (resp. $n$) be the dimension of $D$  (resp. of $G$).
For every $a\in A^{\ast}$ (resp. $b\in B^{\ast}$), fix an
open neighborhood $U_{a}$ of  $a$  (resp.  $V_b$ of $b$) such that $U_{a} $
(resp. $V_b$) is  biholomorphic to a bounded domain in $\C^{m}$  (resp. in
$\C^{n}$).
 For any $0<\delta\leq\frac{1}{2},$  define
\begin{equation} \label{eq_prop3_4}
\begin{split}
U_{a,\delta}&:=\left\lbrace z\in U_{a}:\ \widetilde{\omega}(z, A\cap U_a,   U_a)<\delta  \right\rbrace,\qquad
a\in A\cap A^{\ast},\\
V_{b,\delta}&:=\left\lbrace w\in V_{b}:\  \widetilde{\omega}(w, B\cap V_b,   V_b)<\delta  \right\rbrace,\qquad
b\in B\cap B^{\ast},\\
A_{\delta}&:=\bigcup\limits_{a\in A\cap A^{\ast}} U_{a,\delta},\qquad
B_{\delta}:=\bigcup\limits_{b\in B\cap B^{\ast}} V_{b,\delta},\\
 D_{\delta}&:=\left\lbrace z\in D:\
  \widetilde{\omega}(z,A,D)<1-\delta\right\rbrace,\quad
  G_{\delta}:=\left\lbrace w\in G:\  \widetilde{\omega}(w,B,G)<1-\delta\right\rbrace.
\end{split}
\end{equation}
Observe that
 $U_{a,\delta}$ (resp. $V_{b,\delta}$) is an open neighborhood of $a$  (resp. $b$).
Moreover, one has  the following  inclusion  (which will be implicitly used  in the sequel):
\begin{equation*}
 \X(A\cap A^{\ast},B\cap B^{\ast};D,G)\subset W\cap\widehat{ \widetilde{W}}.
\end{equation*}

\noindent {\bf Step 3:} {\it  The case
where
  $G$ is a bounded open subset  in $\C^n.$}

\smallskip

\noindent {\it Sketchy proof of  Step 3.}
We only describe  the construction of $\hat{f}.$
   For each $a\in A\cap A^{\ast},$  let $f_a:=f|_{\X\left( A\cap U_{a} ,B;  U_{a},G\right)}.$
Since $f\in \mathcal{O}_s(W^{\text{o}},Z),$ we deduce that    $f_a\in
\mathcal{O}_s\Big(\X\left( A\cap U_{a} ,B;  U_{a},G\right),Z\Big).$ Recall that
$U_a$ (resp. $G$) is biholomorphic to a bounded open set in $\C^{m}$ (resp. in $\C^n$). Consequently,  applying Theorem
1
to $f_a$ yields that there is a unique mapping
   $\hat{f}_{a} \in   \mathcal{O}\Big(\widehat{\widetilde{\X}}
 \left(A\cap U_{a} ,B; U_{a},G\right),Z\Big) $
such that
\begin{equation}\label{eq1_prop3}
\hat{f}_{a}(z,w)=f_a(z,w)=f(z,w),\  (z,w)\in
 \X\left(  A\cap A^{\ast}\cap  U_a,B\cap B^{\ast};  U_{a},G\right) .
\end{equation}
Let $0<\delta\leq \frac{1}{2}.$
In virtue of (\ref{eq_prop3_4})--(\ref{eq1_prop3}), we are able to ``glue" the
family $\left(\hat{f}_{a}|_{U_{a,\delta}\times G_{\delta}} \right)_{a\in A\cap A^{\ast}}.$
Let
\begin{equation}\label{eq2_prop3}
\tilde{\tilde{f}}_{\delta}\in \mathcal{O}(A_{\delta}\times
G_{\delta},Z)
\end{equation}
denote the resulting  mapping after the gluing process.
 In virtue of  (\ref{eq1_prop3})--(\ref{eq2_prop3}),
 we are able to define a new mapping  $\tilde{f}_{\delta}$ on $\X\left(A_{\delta},B\cap B^{\ast};D,G_{\delta}
 \right)$ as follows
\begin{equation*}
 \tilde{f}_{\delta}:=
\begin{cases}
 \tilde{\tilde{f}}_{\delta},
  & \qquad\text{on}\  A_{\delta}\times G_{\delta}, \\
  f, &   \qquad\text{on}\ D\times (B\cap B^{\ast})        .
\end{cases}
\end{equation*}
Using this  and  (\ref{eq1_prop3})--(\ref{eq2_prop3}) again, we  see that  $ \tilde{f}_{\delta}\in \mathcal{O}_s\Big(
\X\left(A_{\delta}, B\cap B^{\ast};D,G_{\delta}
 \right),Z\Big),$ and
\begin{equation*}
  \tilde{f}_{\delta}=f\qquad\text{on}\ \X(A\cap A^{\ast},B\cap
 B^{\ast};D,G_{\delta}).
\end{equation*}
 Since $A_{\delta}$ is an open  subset of the complex manifold $D$  and $G_{\delta}$ is biholomorphic to a bounded  open set in $\C^n,$  we are able to  apply  Step 1
to  $ \tilde{f}_{\delta}$ in order to obtain a mapping
  $\hat{f}_{\delta}\in \mathcal{O}\Big(  \widehat{\widetilde{\X}}\left(A_{\delta}, B\cap B^{\ast};D,G_{\delta}
 \right),Z \Big)$  such that
\begin{equation*} 
  \hat{f}_{\delta}=  \tilde{f}_{\delta}\qquad\text{on}\ \X\left(A_{\delta},
  B\cap B^{\ast};D,G_{\delta}\right).
\end{equation*}

We are now in the position to define the desired extension  mapping $\hat{f}.$
Indeed, one  glues
$\left(\hat{f}_{\delta}\right)_{0<\delta\leq\frac{1}{2}}$ together to obtain
$\hat{f}$ in the following way
\begin{equation*} 
\hat{f}:=\lim\limits_{\delta\to 0} \hat{f}_{\delta}\qquad \text{on}\
 \widehat{\widetilde{W}}  .
\end{equation*}
In fact, the equality $\widehat{\widetilde{W}}=\bigcup\limits_{0<\delta<\frac{1}{2}}\widehat{\widetilde{\X}}\left(A_{\delta},
  B\cap B^{\ast};D,G_{\delta}\right)$  follows essentially from  (\ref{eq_prop3_4}).
\hfill $\square$

\smallskip

\noindent {\bf Step 4:} {\it
Completion of the proof of Theorem 5.}

\smallskip

\noindent {\it Sketchy proof of  Step 4.}
 For each $a\in A\cap A^{\ast},$  let $f_a:=f|_{\X\left( A\cap U_{a} ,B;  U_{a},G\right)}.$
Since $f\in \mathcal{O}_s(W^{\text{o}},Z),$ we deduce that    $f_a\in
\mathcal{O}_s\Big(\X\left( A\cap U_{a} ,B;  U_{a},G\right),Z\Big).$ Since $U_a$ is
biholomorphic to a bounded domain in $\C^m,$ we are able to apply  Step 3
 to $f_a.$  Consequently,  there is a  mapping
   $\hat{f}_{a} \in   \mathcal{O}\Big(\widehat{\widetilde{\X}}
 \left(A\cap U_{a} ,B; U_{a},G\right),Z\Big) $
such that
\begin{equation}\label{eq1_prop4}
\hat{f}_{a}(z,w)=f(z,w),\qquad  (z,w)\in
 \X\left(  A\cap A^{\ast}\cap U_a,B\cap B^{\ast};  U_{a},G\right) .
\end{equation}
 Let $0<\delta\leq\frac{1}{2}.$ In virtue of  (\ref{eq1_prop4}), we can ``glue"  the family
 $\left(\hat{f}_{a}|_{U_{a,\delta}\times G_{\delta}} \right)_{a\in A\cap A^{\ast}}$
in order to obtain the resulting   mapping $\tilde{f}^{'}_{\delta}\in \mathcal{O}(A_{\delta}\times
G_{\delta},Z).$

Similarly,  for each $b\in B\cap B^{\ast},$   one obtains a  mapping
 $\hat{f}_{b} \in   \mathcal{O}\Big(\widehat{\widetilde{\X}}
 \left(A, B\cap V_{b} ; D,V_{b}\right),Z\Big) $
such that
\begin{equation}\label{eq2_prop4}
\hat{f}_{b}(z,w) =f(z,w),\qquad (z,w)\in
 \X\left(  A\cap A^{\ast},B\cap B^{\ast}\cap V_b;  D,V_{b}\right) .
\end{equation}
Moreover, one can ``glue"  the family
 $\left(\hat{f}_{b}|_{D_{\delta}\times V_{b,\delta}} \right)_{b\in B\cap B^{\ast}}$
in order to obtain the  resulting  mapping $\tilde{f}^{''}_{\delta}\in \mathcal{O}(D_{\delta}\times
B_{\delta},Z).$

Next, using   (\ref{eq1_prop4})--(\ref{eq2_prop4}) and (\ref{eq_prop3_4}) we can  prove that
\begin{equation*}  
 \tilde{f}^{'}_{\delta}=\tilde{f}^{''}_{\delta}\qquad \text{on}\  A_{\delta}\times
B_{\delta}.
\end{equation*}
Using this
 we are able to define a new  mapping $\tilde{f}_{\delta}:\
 \X\left(A_{\delta}, B_{\delta};D_{\delta},
 G_{\delta}\right)\longrightarrow Z$ as follows
\begin{equation*} 
 \tilde{f}_{\delta}:=
\begin{cases}
\tilde{f}^{'}_{\delta},
  & \qquad\text{on}\  A_{\delta}\times G_{\delta}, \\
  \tilde{f}^{''}_{\delta}, &   \qquad\text{on}\ D_{\delta}\times B_{\delta}        .
\end{cases}
\end{equation*}
Using this formula 
 it can be readily checked that
$\tilde{f}_{\delta}\in \mathcal{O}_s\Big(\X\left(A_{\delta}, B_{\delta};D_{\delta},G_{\delta}
 \right),Z\Big).$ Since  we know from  (\ref{eq_prop3_4}) that  $A_{\delta}$
 (resp.  $B_{\delta}$) is an open subset of $D_{\delta}$  (resp.
 $G_{\delta}$), we are able to apply  Step 2
to  $ \tilde{f}_{\delta}$  for every $0<\delta\leq\frac{1}{2}.$ Consequently,
one
  obtains a   mapping $\hat{f}_{\delta}\in \mathcal{O}\Big(
 \widehat{\X}\left(A_{\delta}, B_{\delta};D_{\delta},G_{\delta}
 \right),Z \Big)$ such that
\begin{equation*}  
  \hat{f}_{\delta}=  \tilde{f}_{\delta}\qquad\text{on}\ \X\left(A_{\delta},
  B_{\delta};D_{\delta},G_{\delta}\right).
\end{equation*}
We are now in the position to define the desired extension  mapping $\hat{f}.$
 \begin{equation*}
\hat{f}:=\lim\limits_{\delta\to 0} \hat{f}_{\delta}\qquad \text{on}\
 \widehat{\widetilde{W}}  .
\end{equation*}
In fact, the equality $\widehat{\widetilde{W}}=\bigcup\limits_{0<\delta<\frac{1}{2}}\widehat{\X}\left(A_{\delta},
  B_{\delta};D_{\delta},G_{\delta}\right)$  follows essentially from  (\ref{eq_prop3_4}).
\hfill $\square$

\section{PROBLEM 1 for the case $A\subset \partial D,$ $B\subset \partial G$}
\label{Section_PROBLEM_1_second_case}

In this  section we present  two particular cases of PROBLEM 1  using two  different systems of approach regions
defined in Subsection  \ref{Subsection_Examples}. These results are obtained in collaboration with Pflug (see \cite{pn1,pn2,pn3}).
Firstly, we start with the case of dimension $1.$

\subsection{System of angular approach  regions}
Our main purpose is to establish a   boundary cross theorem  which is  the optimal version of Theorem 3.
This  constitutes  the first   step of our strategy  to extend  the theory of separately holomorphic mappings.
We will use  the  terminology and the   notation
 in Paragraph  3 of Subsection \ref{Subsection_Examples}.
  More  precisely, if $D$ is  an open set of a  Riemann surface  such that  $D$ is  good on a  nonempty part of $\partial D,$ we equip $D$
    with the  system of angular  approach  regions
 supported on this  part. Moreover, the     notions such as set of positive  length,
 set of zero length, locally pluriregular point 
 which exist on $\partial E$ can be  transferred  to
 $\partial D$   using   conformal mappings  in a local  way  (see \cite{pn2} for  more details).

  \renewcommand{\thethmspec}{Theorem 6  (\cite{pn2})}
  \begin{thmspec}
  Let  $X,\ Y$  be  Riemann  surfaces
  and $D\subset X,$ $ G\subset Y$ open subsets and
  $A$ (resp. $B$) a subset of   $\partial D$ (resp. and more
  $\partial G$)  such that
      $D$ (resp. $G$) is  good on $A$ (resp. $B$) and that both
  $A$ and $B$ are of   positive length.
 Define
 \begin{eqnarray*}
 W&:=&\X(A,B;D,G),\qquad W^{'}:=\X(A^{'},B^{'};D,G),\\
  \widehat{W}  &:= &\left\lbrace  (z,w)\in D\times G:\  \omega(z,A,D)+\omega(w,B,G)<1     \right\rbrace,\\
   \widehat{W^{'}}  &:= &\left\lbrace  (z,w)\in D\times G:\  \omega(z,A^{'},D)+\omega(w,B^{'},G)<1     \right\rbrace,
 \end{eqnarray*}
 where  $A^{'}$  (resp. $B^{'}$) is the set of   points  at  which  $A$ (resp. $B$) is locally pluriregular
 with respect to the system of angular  approach  regions  supported on $A$  (resp.  $B$), and  $\omega(\cdot,A,D),$
 $\omega(\cdot,A^{'},D)$   (resp.  $\omega(\cdot,B,G),$
 $\omega(\cdot,B^{'},G)$) are calculated using the canonical  system of approach regions.

 Then for every function  $f:\ W\longrightarrow \C$ which satisfies the  following conditions:
\begin{itemize}
\item[ (i)]   for  every $a\in A$ the function $f(a,\cdot)|_G$ is  holomorphic and has the angular limit $f(a,b)$
at all points $b\in B,$  and   for  every $b\in B$ the function  $f(\cdot,b)|_D$ is  holomorphic and has the angular limit $f(a,b)$
at all points $a\in A;$
\item[(ii)] $f$ is locally  bounded;
\item[(iii)]    $f|_{A\times B}$ is continuous,
\end{itemize}
there  exists a  unique  function
$\hat{f}\in\mathcal{O}(\widehat{W^{'}},\C)$
which  admits the angular limit $f$ at all points of $W\cap W^{'}.$

If $A$ and $B$  are  Borel  sets  or if $X=Y=\C$ then
 $ \widehat{W} = \widehat{W^{'}}.$
\end{thmspec}

Theorem 6 is  the ``measurable" version of  Theorem 3. Indeed,  the hypotheses of the latter theorem such as
open boundary sets $A$ and $B,$ etc are  now  replaced  by
measurable  boundary sets $A$ and $B,$ etc  in the former theorem.
The  question of optimality  of Theorem 6 has been settled  down  in  \cite{pn4}.

\smallskip

 Our method consists of two steps.
In the first step we  suppose that  $D$ and $G$ are    Jordan domains in $\C.$
In the second one we treat the general case. Now we give a brief outline of the proof.

For every $0<\delta <1$  the set
$D_{\delta}:=\left\lbrace z\in D:\ \omega(z,A,D)<1-\delta \right\rbrace$
(resp.  $G_{\delta}:=\left\lbrace w\in G:\ \omega(w,B,G)<1-\delta \right\rbrace              $)
is  called {\it  a level set} (of the harmonic  measure  $ \omega(\cdot,A,D)$  (resp.   $ \omega(\cdot,B,G)$).
In the first  step, we improve     Gonchar's method  \cite{go1,go2} by   making intensive  use of  Carleman's formula  (see \cite{az}) and of
geometric properties of the level sets of harmonic measures.
 More precisely, by adapting   Gonchar's method  to our  ``measurable" situation, we
meet  some difficulty   concerning the geometry of  $D_{\delta}$ et $G_{\delta}$ which  is  very complicated.
In order to overcome this  situation,  we construct  Jordan domains  with rectifiable boundary which are  contained in  $D_{\delta}$ and $G_{\delta}$
and which touch  the boundary of these level sets   on a set of positive length.
Consequently, the analysis  on the complicated  open sets $D_{\delta}$ and $G_{\delta}$ can be reduced to that  on  certain Jordan domains.

The main ingredient  for the second step is a   mixed cross type theorem.  The idea is to adapt Theorem 1 in the following ``mixed" situation:

{\it $D$ (resp. $G$) is an open set of  a Riemann surface,  $A$ is  an open subset of $D,$ but
$B$ is  a subset  of $ \partial G$ such that  $G$ is  good on $B.$  This  situation explains  the  terminology  ``mixed cross".}

Our key observation is  that    the classical  method of doubly orthogonal bases of Bergman type
that  we  discussed in Section 2  still applies  in the present  mixed  context.
We also  use a recent work of Zeriahi  (see \cite{zer}).

In the second step we  apply this  mixed cross type theorem  in order  to
prove Theorem 6  with $D$  (resp. $G$) replaced by $D_{\delta}$  (resp.  $G_{\delta}$). Then we construct
the solution  for the  original open sets $D$ and $G$ by means of a gluing procedure.
The  method  for the second step (which is called ``the method of level sets")  has appeared for the first time  in \cite{pn1}.
We will discuss it in the next subsection.

\subsection{Canonical  system of  approach  regions}

For every open subset $U\subset \R^{2n-1}$ and every continuous function $h:\  U\longrightarrow\R,$
the graph
$$\left\lbrace z=(z^{'},z_n)=(z^{'},x_n+iy_n)\in \C^n:\ (z^{'},x_n)\in U\quad\text{and}\quad y_n=h(z^{'},x_n)     \right\rbrace$$
is called a {\it topological hypersurface in $\C^n$.}

 Let  $X$ be a complex manifold of dimension $n.$
 A subset $A\subset X$ is said to be  a  {\it topological hypersurface} if,
 for every point $a\in A,$ there is a    local chart  $(U,\phi:\ U\rightarrow\C^n)$ around $a$
 such that  $\phi(A\cap U)$ is a topological hypersurface in $\C^n$

Now let    $D\subset X$  be an open subset and let  $A\subset\partial D$
be an open subset
(with respect to the topology induced on $\partial D$). Suppose in addition that
$A$ is a topological hypersurface.
A point  $a\in A$ is said to be {\it   of type 1  (with respect to $D$)}
if, for every neighborhood $U$ of $a$ there is an open  neighborhood $V$ of $a$ such that $V\subset U$ and  $V\cap D$ is
a domain.  Otherwise, $a$ is said to be {\it   of type 2}.
We see easily that if $a$ is of type 2, then     for every neighborhood $U$ of $a,$
   there are an open neighborhood $V$ of $a$ and
two  domains $V_1,$ $V_2$  such that $V\subset U,$    $V\cap D=V_1\cup V_2$
 and all points in $A\cap V$ are of type 1 with respect to $V_1$ and $V_2.$

 In virtue of  Proposition 3.7 in  \cite{pn3} we  have the following
 \begin{prop}\label{prop_application2_1}
    Let  $X$ be a complex manifold and
     $D$   an open subset of $X.$
  $D$ is  equipped  with  the canonical  system  of  approach  regions.
  Suppose that  $A\subset\partial D$  is an open boundary subset
which is also a topological hypersurface.
Then $A$ is locally  pluriregular and  $\widetilde{A}=A.$
\end{prop}

The main result of this  subsection is

\renewcommand{\thethmspec}{Theorem 7   (\cite{pn3})}
  \begin{thmspec}
Let $X,\ Y$  be two complex manifolds,
  and $D\subset X,$ $ G\subset Y$  two nonempty open sets. $D$ (resp.  $G$) is  equipped  with
  the canonical  system of  approach  regions.
   Let
  $A$ (resp. $B$) be a nonempty open subset of  $\partial D$ (resp.
  $\partial G$) which is also a
   topological hypersurface.
   Define
   \begin{eqnarray*}
   W  &:= &\X(A,B;D,G),\\
    \widehat{W}  &:= &\left\lbrace  (z,w)\in D\times G:\  \omega(z,A,D)+\omega(w,B,G)<1     \right\rbrace.
   \end{eqnarray*}
 Let  $f:\ W\longrightarrow \C$ be   such that:
\begin{itemize}
\item[ (i)]  $f\in\mathcal{C}_s(W,\C)\cap \mathcal{O}_s(W^{\text{o}},\C);$
\item[(ii)]   $f$ is locally bounded on $W;$
\item[ (iii)]   $f|_{A\times B}$ is continuous.
\end{itemize}

Then  there exists a unique function
$\hat{f}\in\mathcal{O}(\widehat{W},\C)$
such that
\begin{equation*}
\lim\limits_{\widehat{W}\ni(z,w)\to(\zeta,\eta)}\hat{f}(z,w)=f(\zeta,\eta),\qquad (\zeta,\eta)\in W.
\end{equation*}
\end{thmspec}

A  weaker version of Theorem 7 where  $D$  (resp. $G$)  is  pseudoconvex open subset of $\C^m$ (resp.  $\C^n$) was previously proved  in \cite{pn1}.
In order to tackle ``arbitrary" complex manifolds  we  follow  our  new approach introduced in  Section  \ref{section_new_approach}  and  \ref{Section_PROBLEM_1_first_case}. The next key technique is to
 apply a   mixed cross type theorem in the following context.

{\it $D$ is  an open  subset of $\C^m$ and  $G$ is  the open unit disc in $\C,$  $A$ is  an open subset of $D$ but
$B$ is  an open connected subset (an arc) of $ \partial G.$}

 The last key technique is to use {\it level sets} of the plurisubharmonic measure (see \cite{pn1,pn2}). More precisely,
we exhaust $D$  (resp.  $G$) by the  level sets of the plurisubharmonic measure
$\omega(\cdot,A,D)$  (resp. $\omega(\cdot,B,G)$),  that is,  by
$D_{\delta}:=\left\lbrace z\in D:\ \omega(z,A,D)<1-\delta \right\rbrace$
(resp.  $G_{\delta}:=\left\lbrace w\in G:\ \omega(w,B,G)<1-\delta \right\rbrace              $)  for  $0<\delta<1.$


 Our method consists of three steps.
In the first step we  suppose that  $G$ is a  domain in $\C^m$ and   $A$ is an open subset of $D.$
In the second step we treat the case  where the pairs  $(D,A)$ and $(G,B)$ are   ``good" enough in the sense of the slicing method.
In the last one we consider the general case. For the first step we  combine the  above mentioned mixed cross  theorem with
the technique of holomorphic discs.  For the second step
 one  applies the slicing method and Theorem 3 \footnote{It is worthy to remark here that  a  weaker version of Theorem 3 will suffice for this argument. Namely, we only need
 Theorem 3 for the case  where $A$ and $B$ are arcs. This  weaker version of Theorem 3 is  also  known under the name
 Dru\.{z}kowski's Theorem (see \cite{dr}). In fact,  we also obtain, by this   way,  a new proof of  Theorem 3 starting
from  Dru\.{z}kowski's Theorem.}.
 The general philosophy is to
prove Theorem   7 with $D$  (resp. $G$) replaced by $D_{\delta}$  (resp.  $G_{\delta}$). Then we construct
the solution  for the  original open sets $D$ and $G$ by means of a gluing procedure (that is, the method of level sets).
In the last step we
transfer the holomorphicity from local situations  to  the global context using Poletsky theory of discs and  Rosay Theorem.

\section{PROBLEM 1 in the general case}
\label{Section_PROBLEM_1_third_case}
  In Section
\ref{Section_PROBLEM_1_first_case} and \ref{Section_PROBLEM_1_second_case}  we have   solved  PROBLEM 1   in some  particular but important cases.
These results  make us  hope that  a  reasonable  solution  to PROBLEM 1 in the  general case may exist. The main purpose of this section
is  to confirm this  speculation. In our work \cite{nv2}  we have introduced  the formulations given in  Section \ref{sec_new_formulations} above
and   developed a unified approach  which  improves the one  given  in Section  \ref{section_new_approach}.
We keep the notation  introduced in Section \ref{sec_new_formulations}, and state the   main results.

\renewcommand{\thethmspec}{Theorem 8  (\cite{nv2})}
  \begin{thmspec}
  Let $X,\ Y$  be two complex manifolds,
  let $D\subset X,$ $ G\subset Y$ be two   open sets, let
  $A$ (resp. $B$) be a subset of  $\overline{ D}$ (resp.
  $\overline{ G}$).  $D$  (resp.  $G$) is equipped with a
  system of approach regions
  $\big(\mathcal{A}_{\alpha}(\zeta)\big)_{\zeta\in\overline{ D},\  \alpha\in I_{\zeta}}$
  (resp.  $\big(\mathcal{A}_{\beta}(\eta)\big)_{\eta\in\overline{ G},\  \beta\in I_{\eta}}$).
Suppose  in addition that  $\widetilde{\omega}(\cdot,A,D)<1$ on $D$ and $\widetilde{\omega}(\cdot,B,G)<1$ on $G.$
   Let $Z$ be a complex analytic  space possessing the  Hartogs extension property.
   Then,
   for every   mapping   $f:\ W\longrightarrow Z$
  which satisfies the following conditions:
   \begin{itemize}
   \item[$\bullet$]    $f\in\mathcal{C}_s(W,Z)\cap \mathcal{O}_s(W^{\text{o}},Z);$
    \item[$\bullet$] $f$ is locally bounded along $\X\big (A\cap\partial D,B\cap\partial G;D,G\big );$ \footnote{
    It follows from  Subsection \ref{subsection_cross}  that $$ \X\big (A\cap\partial D,B\cap\partial G;D,G\big )=
   \big( (D\cup A)\times (B\cap\partial G) \big) \bigcup\big((A\cap\partial D)\times ( G\cup B)\big)
    .$$}
    \item[$\bullet$]          $f|_{A\times B}$ is continuous at all points of
   $(A\cap\partial D)\times (B\cap \partial G),$
    \end{itemize}
     there exists a unique mapping  
$\hat{f}\in\mathcal{O}(\widehat{\widetilde{W}},Z)$ which
 admits  $\mathcal{A}$-limit $f(\zeta,\eta)$ at every point
  $(\zeta,\eta)\in  W\cap   \widetilde{W}.$        
\end{thmspec}

Theorem 8 has an important corollary. Before stating this, we need to
introduce a terminology. A complex manifold $\mathcal{M}$ is said to be
 a {\it Liouville manifold} if $\mathcal{PSH}(\mathcal{M})$ does not
 contain any non-constant bounded above functions. We see clearly that the class
 of Liouville manifolds contains the class of connected compact manifolds.

\renewcommand{\thethmspec}{Corollary 1}
  \begin{thmspec}
We keep the hypotheses and the notation in Theorem 8.
   Suppose in addition that  $G $ is  a Liouville manifold.
Then,
    for every   mapping   $f:\ W\longrightarrow Z$
  which satisfies the following conditions:
   \begin{itemize}
   \item[$\bullet$]    $f\in\mathcal{C}_s(W,Z)\cap \mathcal{O}_s(W^{\text{o}},Z);$
    \item[$\bullet$] $f$ is locally bounded along  $\X\big (A\cap\partial D,B\cap\partial G;D,G\big );$
    \item[$\bullet$]          $f|_{A\times B}$ is continuous at all points of
   $(A\cap\partial D)\times (B\cap \partial G),$
    \end{itemize}
 there is a unique mapping
$\hat{f}\in\mathcal{O}(D\times G,Z)$
which  admits  $\mathcal{A}$-limit $f(\zeta,\eta)$ at every point
  $(\zeta,\eta)\in     W\cap   \widetilde{W}.$        
\end{thmspec}

Corollary 1 follows immediately from Theorem 8 since
$\widetilde{\omega}(\cdot,B,G)\equiv 0.$ This theorem  generalizes, in some sense,  all  results
obtained in  Section
\ref{Section_PROBLEM_1_first_case} and \ref{Section_PROBLEM_1_second_case}.  On  the other hand, we will see  many other applications
of Theorem 8 in Section  \ref{section_application}. We will explain  our unified approach  and  techniques for the proof of Theorem 8 in
the following special ``local" case.

 \begin{prop}  \label{thm_local_euclidean}
 Let $D\subset \C^n,$ $ G\subset \C^m$ be bounded connected open sets.
  $D$  (resp.  $G$) is equipped with a
  system of approach regions
  $\big(\mathcal{A}_{\alpha}(\zeta)\big)_{\zeta\in\overline{D},\  \alpha\in I_{\zeta}}$
  (resp.  $\big(\mathcal{A}_{\alpha}(\eta)\big)_{\eta\in\overline{G},\  \alpha\in I_{\eta}}$).
 Let
  $A$ (resp. $B$) be a nonempty subset of  $\overline{ D}$ (resp.
  $\overline{ G}$) such that  $A$ and $B$ are locally pluriregular.
  Put
    \begin{equation*}
    \begin{split}
    W&:= \X(A,B;D,G),\qquad \overline{W}:=\X(\overline{A},\overline{B};D,G), \\
    \overline{W}^{\text{o}}&:=\X^{\text{o}}(\overline{A},\overline{B};D,G),\qquad
     \widehat{W}:=\widehat{\X}(A,B;D,G).
     \end{split}
    \end{equation*}
   Then,
   for every bounded function  $f:\ \overline{W}\longrightarrow\C$ such that
      $f\in\mathcal{C}_s(\overline{W},\C)\cap \mathcal{O}_s(\overline{W}^{\text{o}},\C)$
      and that $f|_{\overline{A}\times \overline{B}}$ is continuous at all points of $(\overline{A}\cap\partial D)
      \times (\overline{B}\cap\partial G),$
     there exists a unique bounded function
$\hat{f}\in\mathcal{O}(\widehat{W},\C)$ which
  admits  $\mathcal{A}$-limit $f(\zeta,\eta)$ at all points
  $(\zeta,\eta)\in    W.$
\end{prop}

This   result constitutes the core  of the proof of Theorem 8. Indeed,  the latter theorem is, in some sense,
the  ``global " version   of  Proposition  \ref{thm_local_euclidean}. By using   the approach developed in Section \ref{section_new_approach}.
we can  go from  local extensions  to  global ones. In addition, the formulation of Proposition  \ref{thm_local_euclidean} gives  rise to
 Definition  \ref{defi_pluri_measure} of the plurisubharmonic  measure  $\widetilde{\omega}(\cdot,A,D).$
 The core  of  our unified approach will be  presented  below.
Our idea  is to use  an adapted version of Poletsky theory of discs  in order  to reduce Proposition  \ref{thm_local_euclidean}
to the case  where $D$ and $G$ are simply the unit discs and  $A\subset\partial D,$  $B\subset\partial  G$  are measurable sets (that is, a  special case of Theorem   6).

Let us   talk about  the needed version of Poletsky theory of discs.
Let $\mes$ denote the Lebesgue measure on the unit circle $\partial E.$
For a bounded mapping $\phi\in\mathcal{O}(E,\C^n)$  and $\zeta\in \partial E,$
$f(\zeta)$ denotes the angular limit value of $f$ at $\zeta$ if it exists.
A classical theorem of Fatou  says that $\mes\left( \{\zeta\in\partial E:\ \exists f(\zeta)\} \right)=2\pi.$

\begin{prop}\label{Poletsky}
Let $D$ be a bounded  open set in $\C^n,$  $\varnothing\not=A\subset \overline{D},$  $z_0\in D$  and $\epsilon>0.$
Let  $\mathcal{A}$ be a system of approach regions for $D.$
Suppose in addition that  $A$ is locally pluriregular (relative to $\mathcal{A}$) and that
$ \omega(\cdot,A,D)<1$ on $D .$
Then  there exist a  bounded  mapping $ \phi\in\mathcal{O}(E,\C^n)$ and a measurable subset $\Gamma_0\subset \partial E$ with the following  properties:
\begin{itemize}
\item[1)] Every point of $\Gamma_0$ is  a density point of $\Gamma_0,$ $\phi(0)=z_0,$ $\phi(E)\subset  \overline{D},$ $\Gamma_0 \subset \left\lbrace \zeta\in\partial E:\ \phi(\zeta)\in \overline{A}
\right\rbrace,$  and
\begin{equation*}
1-\frac{1}{2\pi}\cdot\mes (\Gamma_0 )<\omega(z_0,A,D)+\epsilon.
\end{equation*}
\item[2)] Let  $f\in\mathcal{C}(D\cup \overline{A},\C)\cap \mathcal{O}(D,\C)$  be such that $f(D)$
 is bounded.  Then there exists a bounded function $g\in\mathcal{O}(E,\C)$ such that
  $g=f\circ \phi$ in  a neighborhood of $0\in E$  and\footnote{ Note  here that  by Part 1),  $(f\circ \phi)(\zeta)$ exists
   for all $\zeta\in \Gamma_0.$}
 $g(\zeta)=(f\circ \phi)(\zeta)$  for all $\zeta\in \Gamma_0.$ Moreover, $g|_{\Gamma_0}\in\mathcal{C}(\Gamma_0,\C).$
 \end{itemize}
\end{prop}

This result is  proved  by adapting the original  discs construction of Poletsky in \cite{po1,po2}.  Recall here that
Poletsky  considered the case  where $A\subset  D$ and  $\mathcal{A}$  is  the canonical  system of approach
regions. But his  method  still   works in our context  by using Montel Theorem on normal families.
It is  worthy to  remark that    $\phi(E)\subset  \overline{D};$ but in general   $\phi(E)\not\subset  D.$

Proposition  \ref{Poletsky} motivates the following
\begin{defi}
\label{candidate}
We keep the hypotheses and notation of Proposition \ref{Poletsky}. Then every pair $(\phi,\Gamma_0)$
satisfying the conclusions 1)--2) of this proposition
 is said to be an {\rm $\epsilon$-candidate for the triplet $(z_0,A,D).$}
\end{defi}

Proposition \ref{Poletsky}  says that  there always exist $\epsilon$-candidates for all triplets $(z,A,D).$
Now  we arrive  at

\smallskip

\noindent{\bf  Sketchy proof of Proposition     \ref{thm_local_euclidean}.}
Firstly, we  give the construction of $\hat{f}.$
Fix  a  point $(z,w)\in\widehat{W},$  we want to determine  the value  $\hat{f}(z,w).$ To do this   let $\epsilon>0$ be such that
 \begin{equation}\label{eq6_Step1_1}
 \omega(z,A,D)+\omega(w,B,G)+2\epsilon<1.
 \end{equation}
   By Proposition \ref{Poletsky} and Definition \ref{candidate}, there is
 an $\epsilon$-candidate $(\phi,\Gamma)$  (resp. $(\psi,\Delta)$)  for $(z,A,D)$  (resp.  $(w,B,G)$).
 Moreover, using the hypotheses, we see that the function $f_{\phi,\psi},$  defined by
 \begin{equation*}
 f_{\phi,\psi}(t,\tau):=f(\phi(t),\psi(\tau)),\qquad  (t,\tau)\in \X\left(\Gamma,\Delta;E,E\right),
 \end{equation*}
satisfies the hypotheses of Theorem 6.
By this theorem, let $\hat{f}_{\phi,\psi}$ be the unique function in $ \widehat{\X}\left(\Gamma,\Delta;E,E\right)$
such that
  \begin{equation*}
   (\Alim \hat{f}_{\phi,\psi})(t,\tau)=f_{\phi,\psi}(t,\tau),\qquad  (t,\tau)\in
  \X^{\text{o}}\left(\Gamma, \Delta;E,E\right),
  \end{equation*}
  where  $\Alim$  is the  angular limit.
  In virtue of (\ref{eq6_Step1_1}) and Proposition \ref{Poletsky}, $(0,0)\in
   \widehat{\X}\left(\Gamma,\Delta;E,E\right).$
   Then  we can define  the value    of the desired extension function $\hat{f}$ at $(z,w)$   as  follows
\begin{equation*}
 \hat{f}(z,w):=\hat{f}_{\phi,\psi}(0,0) .
\end{equation*}

It remains to  prove that the so-defined  $\hat{f}$  possesses the required propertied of Proposition  \ref{thm_local_euclidean}:
namely, $\hat{f}$ is  holomorphic  and   admits the $\mathcal{A}$-limit  $f$ at all points of
  $    W.$

In fact, using the technique of level sets, the holomorphicity of  $\hat{f}$  is reduced  to proving the following mixed  cross version
of Proposition  \ref{thm_local_euclidean}.

\noindent{\bf  Assertion. }{\it  $A$ is  a  measurable subset of  $\partial E$  with $\mes(A)>0,$
    $$D:=\left\lbrace
   w\in E:\ \omega(w,A,E)<1-\delta\right\rbrace   \quad \text{for some}\  \delta:\  0\leq\delta<1,$$  $B$ is  an open subset of an arbitrary complex manifold $G.$}

Using Rosay Theorem, the  case  $\delta=0$ of the  assertion can be  reduced  to  the special case of Theorem 6
where $D$ and $G$ are merely the unit discs and  $A\subset\partial D,$  $B\subset\partial  G$  are measurable sets.

The  case  where $0<\delta<1$ can be  reduced to the previous case by using  conformal mappings
from every connected component of  $D$  onto $E.$  In fact,  all connected components of $D$  are simply connected.
This idea has been developed in  \cite{nv2}, and it is   called
{\it the technique  of conformal mappings}.
The interesting point of  this  proof of the assertion  is  that we avoid  completely the  classical method of
doubly orthogonal bases of Bergman type.

 In order to show that  $\hat{f}$ admits the $\mathcal{A}$-limit  $f$ at all points of
  $    W,$         we make  use  of  an argument  based on Two-Constant Theorem  (see \cite{nv2} for  more details).
\hfill  $\square$

\medskip

In conclusion,   our new
approach illustrates the unified character: ``From local
informations to global extensions". In fact, ``global" results
(i.e.~for general  crosses) can be deduced from ``local"
ones (i.e.~for  boundary crosses defined over the bidisk).

\section{PROBLEM 2} \label{sec_PROBLEM_2}

In the case of  crosses in the interior  context (that is,  $A\subset D$ and $B\subset G$),  one was led to investigate
cross theorems with analytic or pluripolar singularities (see, for
example, \cite{jp2,jp3,jp4,jp5} and the references therein). The
starting point to this kind of questions was the so--called range
problem in the theory of mathematical tomography (for more details
see \cite{ok1}). To be more precise one had to describe the range of
the exponential Radon transform $R_\mu$, $\mu\neq 0$,

$$
\cal C^\infty_c(\R^2,\R)\ni h\overset{R_\mu} \mapsto
\int_{x\cdot\omega=p}h(x)\exp(\mu x\cdot\omega^\bot)d\Lambda_1(x),
$$
where $\omega=(\sin \alpha,\cos\alpha)\in\S^1$, $p\in\R$,
$\omega^\bot=(-\sin \alpha,\cos\alpha)$, and where ``$\cdot$''
means the standard scalar product in $\R^2$  and $d\Lambda_1$ denotes  the one-dimensional Lebesgue measure.

Then the natural question arises whether there also exists a general cross theorem with singularities.
Namely, does there exist  a general version of Theorem 2 in the spirit of Theorem 8?
In other words,  we want to solve  PROBLEM 2   when
$Z$ is a complex analytic  space possessing the  Hartogs extension property.

We have recently  obtained, in collaboration
with P. Pflug (see \cite{pn5,np1}), a  reasonable  solution to the problem.
Our idea is   to follow
the strategy as in the  case without singularities. Namely, we
investigate first the ``local" case where the boundary crosses are defined
over the bidisk, then   we pass from this case  to the global one.

By using  an idea of Jarnicki and Pflug  in  \cite{jp3,jp6},
applying  the technique of conformal mappings (see  the end of Section \ref{Section_PROBLEM_1_third_case}),
using  the technique of level sets
and using the results of
Chirka \cite{Chi}, Imomkulov--Khujamov  \cite{ik}  and Imomkulov \cite{im},
 we obtain the following    ``measurable" version  with  singularities
of Theorem 3.

\renewcommand{\thethmspec}{Theorem 9  (\cite{pn5})}
  \begin{thmspec}
Let  $D=G=E$ and  let $A\subset\partial D, $  $B\subset\partial G$
be  measurable  subsets such that $\mes(A)>0,$
$\mes(B)>0.$  Suppose that $D$ and $G$ are equipped with the system of angular approach regions.
Consider the cross $W:=\X(A,B;D,G).$
Let $M$ be a relatively closed subset of $W$  such that
\begin{itemize}
\item[$\bullet$] $M_a$ is   polar (resp.
discrete)  in $G$ for all $a\in A$ and
 $M^b$
is   polar (resp. discrete)   in $D$  for all $b\in B;$ \footnote{ In other words,
$M$ is  polar (resp. discrete) in fibers over $A$ and $B.$}
\item[$\bullet$] $M\cap (A\times B)=\varnothing.$
\end{itemize}
Then there exists a relatively closed pluripolar  subset (resp. an analytic subset)  $\widehat{M}$ of $\widehat{W}$
with the following two properties:
\begin{itemize}
\item[(i)]
The set of end-points of $\widehat{W}\setminus \widehat{M}$ contains
 $\big((A^{'}\times G)\bigcup(  D\times B^{'})\big)\setminus M,$
where $A^{'}$  (resp.  $B^{'}$) denotes  the set of density points  of $A$  (resp.  of $B$).
\item[(ii)] Let $f:\ W\setminus M \longrightarrow\C$
 be a locally bounded
function  such that
\begin{itemize}
\item[$\bullet$]  for all $a\in A,$
 $f(a,\cdot)|_{G\setminus M_a}$   is
holomorphic and admits the angular limit $f(a,b)$ at all points
$b\in B;$
 \item[$\bullet$] for all $b\in B,$
$f(\cdot,b)|_{D\setminus M^b}$ is holomorphic and admits the
angular limit $f(a,b)$ at all points $a\in A;$
\item[$\bullet$]
$f|_{A\times B}$ is measurable.
\end{itemize}
  Then
there is a unique function
$\hat{f}\in\mathcal{O}(\widehat{W}\setminus \widehat{M},\C)$ such
that $\hat{f}$ admits the angular limit $f$ at  all points of
$\big((A^{''}\times G)\bigcup(  D\times B^{''})\big)\setminus M,$
where $A^{''}$  (resp.  $B^{''}$) is  a  subset of $A^{'}$  (resp.  of $B^{'}$)
with  $\mes(A^{'}\setminus A^{''})=0$  (resp.   $\mes(B^{'}\setminus B^{''})=0$).
\end{itemize}

Moreover,   if $M=\varnothing,$ then $\widehat{M}=\varnothing.$
\end{thmspec}

The itinerary to go from Theorem 9 to its global  version is much harder  than  that in the case  without singularities.
The  difficulty  arises when  we want to  show that  $\hat{f}$ admits the desired $\mathcal{A}$-limit.
  In the case  without singularities
this procedure  works well  because  we can use an argument  based on Two-Constant Theorem.
But this  is  not available  any more in the case  with singularities. In \cite{np1} we  have found
a way to overcome this  difficulty by  using some special mixed cross theorems with singularities.

 Recall that a subset $S$ of a complex  manifold $\mathcal{M}$ is said to be  {\it thin} if  for    every point $x\in\mathcal{M}$  there  are a connected
 neighborhood $U=U(x)\subset\mathcal{M}$ and  a holomorphic  function $f$ on $U,$ not identically  zero,
 such that  $U\cap S\subset f^{-1}(0).$ We are now ready to state our main result.

\renewcommand{\thethmspec}{Theorem 10  (\cite{np1})}
  \begin{thmspec}
 Let $X,\ Y$  be two complex manifolds,
  let $D\subset X,$ $ G\subset Y$ be two open sets, let
  $A$ (resp. $B$) be a subset of  $\overline{ D}$ (resp.
  $\overline{ G}$).  $D$  (resp.  $G$) is equipped with a
 system of approach regions
  $\big(\mathcal{A}_{\alpha}(\zeta)\big)_{\zeta\in\overline{ D},\  \alpha\in I_{\zeta}}$
 (resp.  $\big(\mathcal{A}_{\beta}(\eta)\big)_{\eta\in\overline{ G},\  \beta\in I_{\eta}}$).
Suppose  in addition that  $A=A^{\ast}$ and $B=B^{\ast}$
 \footnote{     It is  worthy to note that this assumption is  not so restrictive since
 we know  from  Subsection \ref{Subsection_approach_regions} that $A\setminus A^{\ast}$  and   $ B\setminus
 B^{\ast}$
 are locally pluripolar for arbitrary sets $A\subset\overline{D},$ $B\subset\overline{G}$.    }
 and that $\widetilde{\omega}(\cdot,A,D)<1$ on $D$ and $\widetilde{\omega}(\cdot,B,G)<1$ on $G.$
   Let $Z$ be a complex analytic  space possessing the  Hartogs extension property.
  Let $M$ be  a  relatively closed subset of  $W$  with the following properties:
\begin{itemize}
\item[$\bullet$] $M$ is  thin  in  fibers  (resp.   locally   pluripolar in fibers) over  $A$ and  over $B;$
\item[$\bullet$] $M\cap \big((A\cap\partial D)\times B\big)=
M\cap \big( A\times (B\cap\partial G)\big)=\varnothing.$
\end{itemize}
Then  there exists a relatively closed    analytic  (resp. a  relatively  closed locally  pluripolar) subset
$\widehat{M}$ of $\widehat{\widetilde{W}}$ such that  $\widehat{M}\cap \widetilde{W}\subset M$ \footnote{ Note that
if $\widetilde A\cap
D=\varnothing$ and $\widetilde B\cap G=\varnothing$, then this intersection is empty.} and that
 $\widetilde{W}\setminus M\subset  \End(\widehat{\widetilde{W}}\setminus \widehat{M})$  and that
 for  every mapping    $f:\ W\setminus M\longrightarrow Z$
  satisfying the following  conditions:
  \begin{itemize}
   \item[(i)]    $f\in\mathcal{C}_s(W\setminus M,Z)\cap \mathcal{O}_s(W^{\text{o}}\setminus M,Z);$
    \item[(ii)] $f$ is locally bounded   along  $\X(A\cap\partial D,B\cap\partial G;D,G)
    \setminus M;$
    \item[(iii)]          $f|_{(A\times B)\setminus M}$ is  continuous at  all  points of
   $(A\cap\partial D)\times (B\cap \partial G),$
    \end{itemize}
     there exists  a unique  mapping
$\hat{f}\in\mathcal{O}(\widehat{\widetilde{W}}\setminus \widehat{M}
,Z)$ which  admits the $\mathcal{A}$-limit $f(\zeta,\eta)$ at  every
point
  $(\zeta,\eta)\in    \widetilde{W}\setminus M  .$
\end{thmspec}

\section{Some applications}\label{section_application}
In \cite{nv2} the  author gives various  applications of
Theorem  8  using three systems of
approach regions. These  are the canonical one,  the system of
angular approach regions and the system of conical approach regions.
We only give here  some applications of Theorem 10 for  the system of conical
approach regions.   We leave  the reader to  treat the two first
cases, that is, to translate Theorem 6 and 7
into the new context of  Theorem  10.

Let  $X$  be  an arbitrary  complexe manifold  and $D\subset  X$  an
open subset.
 We say  that  a set   $A\subset \partial D$ is {\it locally contained in a  generating manifold} if  there  exist  an (at most countable)
 index set $J\not=\varnothing,$ a  family  of open subsets
 $(U_j)_{j\in J}$  of $X$ and  a  family  of
   {\it generating manifolds} \footnote{ A differentiable  submanifold  $\mathcal{M}$ of a complex manifold
  $X$ is  said to be a {\it generating
  manifold}
  if  for all  $\zeta\in\mathcal{M},$   every  complex vector subspace of  $T_{\zeta}X$ containing   $T_{\zeta}\mathcal{M}$
coincides  with
    $T_{\zeta}X.$}   $(\mathcal{M}_j)_{j\in J}$ such that
    $A\cap U_j\subset \mathcal{M}_j,$  $j\in J,$ and that $A\subset  \bigcup_{j\in J}  U_j.$
The dimensions of $\mathcal{M}_j$ may  vary according  to $j\in J.$

Suppose that  $A\subset \partial D$ is  locally contained in a
generating manifold. Then we say that $A$ is {\it of positive size}
if under the above notation  $\sum_{j\in J}\mes_{
\mathcal{M}_j}(A\cap U_j)>0,$ where $\mes_{\mathcal{M}_j}$ denotes
the    Lebesgue measure on $\mathcal{M}_j.$
 A point $a\in A$ is  said  to be  a {\it density point} relative to $A$ if it is a density point relative to $A\cap U_j$ on $\mathcal{M}_j$ for some $j\in J.$
 Denote  by  $A^{'}$ the set of all density points   relative to $A.$

 Suppose now that  $A\subset \partial D$ is  of positive size.
 We  equip $D$  with the  system of conical approach regions   supported on   $A.$
  Using the   works of B. Coupet  and   B. J\"{o}ricke  (see \cite{co,jr}),
 one can show that \footnote{  A complete proof will be  available in \cite{nv3}.}   $A$  is  locally  pluriregular  at all density points
  relative to $A$ and $A^{'}\subset \widetilde{A}.$  Consequently, it  follows from  Definition
\ref{defi_pluri_measure}  that
\begin{equation*}
\widetilde{\omega}(z,A,D)\leq  \omega(z,A^{'},D),\qquad  z\in D.
\end{equation*}
This estimate, combined with  Theorem 10, implies the  following
result.

\renewcommand{\thethmspec}{Corollary 2}
  \begin{thmspec}
  Let $X,\ Y$  be  two complex  manifolds,
  let $D\subset X,$ $ G\subset Y$ be  two connected  open sets, and let
  $A$ (resp. $B$) be  a  subset of  $\partial D$ (resp.
  $\partial G$).
   $D$  (resp.  $G$) is  equipped with a
  system of conical approach  regions
  $\big(\mathcal{A}_{\alpha}(\zeta)\big)_{\zeta\in\overline{ D},\  \alpha\in I_{\zeta}}$
  (resp.  $\big(\mathcal{A}_{\beta}(\eta)\big)_{\eta\in\overline{ G},\  \beta\in I_{\eta}}$)  supported on  $A$  (resp. on $B$).
   Suppose  in  addition that
    $A$  and  $B$  are  of positive size.
   Define
   \begin{eqnarray*}
   W^{'}  &:= &\X(A^{'},B^{'};D,G),\\
    \widehat{W^{'}}  &:= &\left\lbrace  (z,w)\in D\times G:\  \omega(z,A^{'},D)+\omega(w,B^{'},G)<1     \right\rbrace,
   \end{eqnarray*}
     where $A^{'}$  (resp.  $B^{'}$) is  the set of  density  points   relative to $A$  (resp.  $B$).
Let $M$ be  a  relatively  closed subset of  $W$  with the following properties:
\begin{itemize}
\item[$\bullet$] $M$ is   thin in  fibers  (resp.   locally   pluripolar in fibers) over  $A$ and  over $B;$
\item[$\bullet$] $M\cap (A\times B)=\varnothing.$
\end{itemize}
Then  there exists a relatively closed    analytic  (resp. a  relatively  closed
locally  pluripolar) subset
 $\widehat{M}$ of $\widehat{W^{'}}$ such that
for  every mapping    $f:\ W\setminus M\longrightarrow Z$
   satisfying the following  conditions:
   \begin{itemize}
   \item[(i)]    $f\in\mathcal{C}_s(W\setminus M,Z)\cap \mathcal{O}_s(W^{\text{o}}\setminus M,Z);$
    \item[(ii)] $f$ is locally bounded   along  $\X(A,B;D,G)
    \setminus M;$
    \item[(iii)]          $f|_{(A\times B)}$ is  continuous,
    \end{itemize}
     there exists  a unique  mapping
$\hat{f}\in\mathcal{O}(\widehat{W^{'}}\setminus \widehat{M} ,Z)$
which admits the $\mathcal{A}$-limit $f(\zeta,\eta)$ at  every point
  $(\zeta,\eta)\in  (W\cap   W^{'})\setminus M  .$
\end{thmspec}

The second  application is  a very general mixed cross   theorem.

 \renewcommand{\thethmspec}{Corollary 3}
  \begin{thmspec}
  Let $X,\ Y$  be  two complex  manifolds,
  let $D\subset X,$ $ G\subset Y$ be  connected open sets,  let
  $A$ be a subset of  $\partial D,$   and let  $B$  be  a subset of
  $ G$.
   $D$ is  equipped with the
  system of conical approach regions
  $\big(\mathcal{A}_{\alpha}(\zeta)\big)_{\zeta\in\overline{ D},\  \alpha\in I_{\zeta}}$  supported on  $A$
   and  $G$  is equipped with the canonical   system of
   approach  regions
  $\big(\mathcal{A}_{\beta}(\eta)\big)_{\eta\in\overline{ G},\  \beta\in I_{\eta}}$.
   Suppose in addition that
    $A$  is  of positive size and  that $B=B^{\ast}\not=\varnothing$.
   Define
    \begin{eqnarray*}
   W^{'}  &:= &\X(A^{'},B;D,G),\\
    \widehat{W^{'}}  &:= &\left\lbrace  (z,w)\in D\times G:\  \omega(z,A^{'},D)+\omega(w,B,G)<1     \right\rbrace,
   \end{eqnarray*}
    where $A^{'}$  is the set of  density points  relative to $A.$
Let $M$ be  a  relatively  closed subset of  $W$  with the following properties:
\begin{itemize}
\item[$\bullet$] $M$ is  thin  in  fibers  (resp.  locally   pluripolar in fibers) over  $A$ and  over $B;$
\item[$\bullet$] $M\cap (A\times B)=\varnothing.$
\end{itemize}
 Then  there exists a relatively closed  analytic  (resp. a  relatively  closed locally  pluripolar) subset
 $\widehat{M}$ of $\widehat{W^{'}}$ such that $   W^{'}\setminus M\subset \End(\widehat{W^{'}}\setminus
 \widehat{M})$
 and that
for  every mapping    $f:\ W\setminus M\longrightarrow Z$
   satisfying the following  conditions:
  \begin{itemize}
   \item[(i)]    $f\in\mathcal{C}_s(W\setminus M,Z)\cap \mathcal{O}_s(W^{\text{o}}\setminus M,Z);$
    \item[(ii)] $f$ is locally bounded   along  $(A\times G)
    \setminus M,$
    \end{itemize}
     there exists  a unique  mapping
$\hat{f}\in\mathcal{O}(\widehat{W^{'}}\setminus \widehat{M} ,Z)$
which admits the $\mathcal{A}$-limit $f(\zeta,\eta)$ at  every point
  $(\zeta,\eta)\in    W^{'}\setminus M  .$
\end{thmspec}

Recently,  Sadullaev  and  Imomkulov (see \cite{is}) have  obtained some  similar results, but not so general  as    Corollary 3. In fact, they introduced the inner plurisubharmonic measure for boundary sets and formulated their results using this function.

\section{Concluding remarks and open questions} \label{last_section}
We collect here  some  open questions  which  seem to be  of interest for the future developments
of  the theory of separately holomorphic mappings.

\smallskip

\noindent{\bf Question 1.}  {\it  Study the optimality  of Theorem 8 and 10.
}

\smallskip

\smallskip

\noindent{\bf Question 2.}  {\it  Investigate  {\rm PROBLEM 1} and   {\rm  2}  when the ``target  space" $Z$ does not possess the Hartogs extension property.
}

\smallskip

\noindent{\bf Question 3.}  {\it Study {\rm PROBLEM 1} when $D$ and $G$ are not necessarily open subsets of $X$ and $Y.$
Here $\mathcal{O}(D,Z)$ denotes the set of all holomorphic  mappings $f:\ U\rightarrow Z,$ where $U=U_f$ is an
open neighborhood of $D$ in $X$ that  depends on $f.$}

Some  results concerning Question 2 could be found in \cite{iv2,iv3,iv4}.
Question 3 has some relations with Sibony's work in \cite{ne}.

We think that  new  tools and  new  ideas  need to be  introduced  in order to  solve  these  questions.

\end{document}